\documentclass[final]{siamltex}
\usepackage{animate}
\usepackage{amsmath}
\usepackage{amssymb}
\usepackage{url}
\usepackage{graphics}
\usepackage{subfigure}
\usepackage{graphicx}
\usepackage{textcomp}
\usepackage{mathrsfs}
\usepackage{epstopdf}
\usepackage{algorithmic, algorithm}
\usepackage{tabu}

\makeatletter

\def\zapcolorreset{\let\reset@color\relax\ignorespaces}
\def\colorrows#1{\noalign{\aftergroup\zapcolorreset#1}\ignorespaces}

\makeatother

\newtheorem{remark}{Remark}[section]

\title{Superconvergent Two-Grid Methods for Elliptic Eigenvalue Problems}

\author{Hailong Guo\footnotemark[2]
\and Zhimin Zhang\footnotemark[2]
\and  Ren Zhao\footnotemark[2]}

\begin{document}

\maketitle
\renewcommand{\thefootnote}{\fnsymbol{footnote}}
\footnotetext[2]{Department of Mathematics, Wayne State University,
      Detroit, MI 48202 ({\tt guo@math.wayne.edu}, {\tt zzhang@math.wayne.edu},
    {\tt renzhao@wayne.edu}).
    This work is supported in part by the US National Science Foundation
    through grant 1115530.}
    \renewcommand{\thefootnote}{\arabic{footnote}}

\begin{abstract}
Some numerical algorithms for elliptic eigenvalue problems are proposed, analyzed, and numerically tested.
The methods combine advantages of the two-grid algorithm
[J. Xu and A. Zhou, Math. Comp, 70(2001), 17--25],
 two-space method [M.R. Racheva and A. B. Andreev, Comput.
Methods Appl. Math., 2(2002), 171--185], the shifted inverse power method
 [X. Hu and X. Cheng,  Math. Comp., 80(2011), 1287--1301;
  Y. Yang and H. Bi, SIAM J. Numer. Anal, 49(2011), 1602--1624],
and the polynomial preserving recovery technique  [Z. Zhang and
A. Naga, SIAM J. Sci.  Comput., 26(2005), 1192--1213].
Our new algorithms compare favorably with some existing methods and enjoy superconvergence property.
\end{abstract}



\pagestyle{myheadings}
\thispagestyle{plain}

\section{Introduction}
A tremendous variety of science and engineering applications, e.g. the  buckling
of columns and shells and  the vibration of elastic bodies, contain  models of
eigenvalue problems of partial differential equations. A recent survey article \cite{grebenkov2013} of SIAM Review
listed 515 references on theory and application of the Laplacian eigenvalue problem.
As one of the most popular numerical methods,  finite element
method has attracted considerable attention in  numerical solution
of eigenvalue problems. A priori error estimates for the finite element
approximation of eigenvalue problems have been investigated
by many authors,  see e.g.,  Babu$\check{\text{s}}$ka and
Osborn \cite{babuska1989, babuska1991},  Chatelin \cite{chatelin},
Strang and Fix \cite{strang}, and references cited therein.

To reduce the computational cost of eigenvalue problems,
Xu and Zhou introduced a two-grid discretization scheme \cite{xuzhou2001}.
Later on, similar ideas were applied to non self-adjoint eigenvalue problems
\cite{kolman} and semilinear elliptic eigenvalue problems \cite{chien}.
Furthermore, it also has been generalized to three-scale
discretization \cite{gao2008}
and  multilevel discretization \cite{linxie}. Recently, a new
shifted-inverse power method based two-grid scheme was proposed in
\cite{hucheng, yangbi2011}.

To improve accuracy of eigenvalue approximation,  many
methods have been proposed. In \cite{shen2006}, Shen and Zhou introduced
a defect correction scheme based on  averaging recovery, like  a
global $L^2$ projection and a  Cl$\acute{\text{e}}$ment-type operator.
In \cite{nagazhangzhou2006}, Naga, Zhang, and Zhou used Polynomial Preserving
Recovery to enhance eigenvalue approximation.  In \cite{wuzhang2009},
Wu and Zhang further showed polynomial preserving recovery can even
enhance eigenvalue approximation on adaptive  meshes.
The idea was further studied in \cite{mengzhang2012, fang2013}.
Alternatively, Racheva and Andreev proposed a two-space method to achieve
better eigenvalue approximation  \cite{racheva} and it was  also
applied to biharmonic eigenvalue problem \cite{andreev}.

In this paper, we propose some fast and efficient solvers
for elliptic eigenvalue problems. We combine ideas of the
two-grid method,  two-space method,  shifted-inverse power method,
and PPR recovery enhancement to design our new algorithms.
 The first purpose is to introduce two superconvergent
two-grid methods for eigenvalue problems.  Our first algorithm is a combination
of the shifted-inverse power based two-grid  scheme
\cite{hucheng, yangbi2011} and
polynomial preserving recovery enhancing technique \cite{nagazhangzhou2006}.
It is worth to point out that the first algorithm can also be seen as
post-processed  two-grid scheme.  We should mention that \cite{liu2011}
also considered
postprocessed two-scale finite element discretization for elliptic
partial operators including boundary value problems and eigenvalue problems.
However, their method is limited to tensor-product domains. Our method
works for arbitrary domains and hence is more general.
The second algorithm can be viewed as a combination of the two-grid
scheme \cite{hucheng, yangbi2011} and the two-space method \cite{racheva,
andreev}. It  can be thought as a special $hp$ method.
The new proposed methods enjoy all advantages of the above methods :
low computational cost and superconvergence.

Solutions of practical problems are often suffered from low regularity.
Adaptive finite element method(AFEM) is a fundamental tool to overcome such
difficulties. In the context of adaptive finite element method for elliptic
eigenvalue problems, residual type a posteriori error  estimators are analyzed
in \cite{duran2003, hoppe2010, larson, verfurth} and
recovery type a posteriori error  estimators are  investigated by
\cite{mao2006, wuzhang2009, liuyan2012}.
For all adaptive methods mentioned above,
 an algebraic eigenvalue  problem has to be solved during every iteration, which
 is very time consuming. This cost dominates the computational cost of
 AFEM and usually is ignored. To reduce computational cost,
Mehrmann and Miedlar \cite{mehrmann} introduced a new adaptive  method which
only  requires an inexact solution of algebraic eigenvalue equation on each
iteration by only performing a few iterations of Krylov subspace solver.
Recently, Li and Yang \cite{liyang} proposed an adaptive
finite element method based on multi-scale discretization for eigenvalue problems
and Xie \cite{xie} introduced a type of adaptive finite element  method
based on  the multilevel correction scheme. Both methods only solve  an
eigenvalue problem on the coarsest mesh and solve boundary value  problems on
adaptive refined meshes.

The second purpose of this paper is to propose two multilevel
adaptive methods. Using our methods, solving an eigenvalue problem by
AFEM will not be more difficult  than solving  a boundary value problem by AFEM.
The most important feature which  distinguishes them  from the methods
in \cite{liyang, xie} is that superconvergence
of eigenfunction approximation and ultraconvergence (two order higher)
of eigenvalue approximation can be numerically observed.

The rest of this paper is organized as follows. In Section 2, we introduce
finite element discretization of elliptic eigenvalue problem and polynomial
preserving recovery. Section 3 is devoted to  presenting two superconvergent
two-grid methods and their error estimates. In Section 4, we propose two
multilevel adaptive methods. Section 5 gives some numerical examples to
demonstrate efficiency of our new methods and finally some conclusions are
draw in Section 6.

\section{Preliminary} In this section, we first introduce the model eigenvalue
problem and its conforming finite element discretization.  Then, we
give a simple description of polynomial preserving recovery for linear element.
\subsection{ A PDE eigenvalue problem and its finite element discretization} Let $\Omega\subset
\mathbb{R}^2$ be a polygonal domain with Lipschitz continuous boundary
$\partial \Omega$. Throughout this article, we shall use the standard notation
for classical Sobolev spaces and their associated norms, seminorms, and
inner products as in \cite{brenner, ciarlet}. For a subdomain $G$ of $\Omega$,
$W^{k, p}(G)$ denotes the classical Sobolev space with norm
$\|\cdot\|_{k, p, G}$ and the seminorm $|\cdot|_{k, p, G}$. When $p=2$,
$H^{m}(G):=W^{m, 2}(G)$  and the index $p$ is omitted. In this article, the
letter $C$, with or without subscript, denotes a generic constant which is
independent of mesh size $h$ and may not be the same at each occurrence.
To simplify notation, we denote $X\leq CY$ by $X \lesssim Y$.

Consider the following second order self adjoint elliptic eigenvalue problem:
\begin{equation}
    \begin{cases}
     -\nabla \cdot (\mathcal{D} \nabla u) + cu = \lambda u,
     & \forall x \in \Omega,\\
	u|_{\partial \Omega} = 0.& \\
    \end{cases}
    \label{equ:model}
\end{equation}
where $\mathcal{D}$ is a $2\times 2$ symmetric positive definite matrix
and $c \in L^{\infty}(\Omega)$. Define a bilinear form $a(\cdot, \cdot):
H^{1}(\Omega) \times H^{1}(\Omega) \rightarrow \mathbb{R}$ by
\begin{equation*}
    a(u, v) = \int_{\Omega}(\mathcal{D}\nabla u\cdot \nabla v + cuv)dx.
\end{equation*}
Without loss of generality, we may assume that $c\ge 0$. It is easy to see that
\begin{equation*}
    a(u, v) \leq C \|u\|_{1, \Omega}\|v\|_{1, \Omega}, \quad \forall
    u, v \in H^{1}_{0}(\Omega),
\end{equation*}
and
\begin{equation*}
    a(u, u) \geq \alpha \|u\|_{1, \Omega}^2, \quad
    \forall u \in H^{1}_{0}(\Omega).
\end{equation*}
Define $\|\cdot\|_{a, \Omega} = \sqrt{a(\cdot, \cdot)}$. Then
$\|\cdot\|_{a, \Omega}$ and $\|\cdot\|_{1, \Omega}$ are two equivalent norms
in $H^{1}_{0}(\Omega)$.

The variational formulation of \eqref{equ:model} reads as: Find
$(\lambda, u)\in \mathbb{R}\times H^{1}_{0}(\Omega)$ with $u \neq 0$ such that
\begin{equation}
	a(u, v) = \lambda (u, v) , \quad  \forall v \in H^{1}_{0}(\Omega).
    \label{equ:variational}
\end{equation}
It is well known that \eqref{equ:variational} has a countable sequence of
real eigenvalues $0 < \lambda_1 \leq \lambda_2 \leq \lambda_3\leq \cdots
\rightarrow \infty$ and corresponding eigenfunctions $u_1, u_2, u_3, \cdots$
which can be assumed to satisfy $a(u_i, u_j) = \lambda_i (u_i, u_j) =
\delta_{ij}$. In the sequence $\{\lambda_j\}$ , the $\lambda_i$ are repeated
according to geometric multiplicity.

Let $\mathcal{T}_h$ be a conforming triangulation of the domain $\Omega$ into
triangles $T$ with diameter $h_{T}$ less than or equal to $h$.
Furthermore, assume
$\mathcal{T}_h$ is shape regular \cite{ciarlet}.
 Let $r \in \{1, 2\}$ and define the continuous finite element space of order $r$ as
\begin{equation*}
    S^{h, r} = \left\{ v \in C(\bar{\Omega}): v|_{T} \in \mathbb{P}_{r}(T),
    \forall T\in \mathcal{T}_h\right\} \subset H^{1}(\Omega),
\end{equation*}
where $\mathbb{P}_r(T)$ is the space of polynomials of degree less than or
equal to $r$ over $T$. In addition, let
$S^{h, r}_{0} = S^{h, r}  \cap H^{1}_0(\Omega)$. In most cases, we shall
use linear finite element space and hence denote $S^{h, 1}$ and $S^{h, 1}_0$
by $S^h$ and $S^h_0$ to simplify notation. The finite element discretization
of \eqref{equ:model} is : Find
$(\lambda_h, u_h) \in \mathbb{R}\times S^{h, r}_0$ with $u_h \neq 0$ such that
\begin{equation}
	a(u_h, v_h) = \lambda_h (u_h, v_h) , \quad
	\forall v_h \in S^{h, r}_{0}.
    \label{equ:fem}
\end{equation}
Similarly, \eqref{equ:fem} has a finite sequence of eigenvalues
$0 <  \lambda_{1, h} \leq \lambda_{2, h} \leq \cdots \leq \lambda_{n_h, h}$ and
corresponding eigenfunctions $u_{1, h}, u_{2, h}, \cdots, u_{n_h, h}$
which can be chosen to satisfy  $a(u_{i, h}, u_{j, h}) = \lambda_{i, h}
(u_{i, h}, u_{j, h}) = \delta_{ij}$ with $i, j = 1, 2, \cdots, n_h$ and
$n_h = \dim S^{h, r}_{0}$.


Suppose that the algebraic multicity of $\lambda_i$ is equal to $q$, i.e.
$\lambda_{i} = \lambda_{i+1} = \cdots  = \lambda_{i+q-1}$.
Let $M(\lambda_i)$ be the space spanned by all eigenfunctions corresponding
to $\lambda_i$. Also, let $M_h(\lambda_h)$ be the direct sum of
eigenspaces corresponding to all eigenvalue $\lambda_{i,h}$ that convergences
to $\lambda_i$.


For the above conforming finite element discretization, the following
 result has been established by many authors
\cite{babuska1991, xuzhou2001, yangbi2011}.

\begin{theorem}
    Suppose $M(\lambda_{i})\subset  H^1_0(\Omega)\cap H^{r+1}(\Omega)$.
    Let $\lambda_{i, h}$  and $\lambda_i$ be the $i$th eigenvalue of
    \eqref{equ:fem} and  \eqref{equ:variational}, respectively.
    Then
    \begin{equation}
	\lambda_{i} \le \lambda_{i, h} \le \lambda_{i} + C h^{2r}.
	\label{equ:eigenvalueerror}
    \end{equation}
    For any eigenfunction $u_{i, h}$ corresponding to $\lambda_{i, h}$
    satisfying $\|u_{i, h}\|_{a, \Omega} = 1$, there exists
    $u_{i} \in M(\lambda_i)$ such that
    \begin{equation}
	\|u_{i}-u_{i, h}\|_{a, \Omega} \le Ch^r.
	\label{equ:eigenfunctionerror}
    \end{equation}
    \label{thm:convergence}
\end{theorem}

Before ending this subsection, we present an important  identity
\cite{babuska1991}  of eigenvalue and  eigenfunction approximation.
\begin{lemma}
    Let $(\lambda, u)$ be the solution of \eqref{equ:variational}. Then
    for any $w \in H^1_0(\Omega)\backslash \{0\}$, there holds
    \begin{equation}
	\frac{a(w, w)}{(w, w)} - \lambda =
	\frac{a(w-u, w-u)}{(w, w)} - \lambda
	\frac{(w-u, w-u)}{(w, w)}.
	\label{equ:identity}
    \end{equation}
    \label{lem:identity}
\end{lemma}
This identity will play an important role in our superconvergence analysis.

\subsection{Polynomial Preserving Recovery} Polynomial Preserving
Recovery (PPR) \cite{zhangnaga2005, nagazhang2004, nagazhang2005}
is an important  alternative of the famous Superconvergent Patch
Recovery proposed by  Zienkiewicz and Zhu \cite{zhu1992}.
Let $G_h: S^h \rightarrow S^h\times S^h$ be the PPR operator and $u_h$ be
a function in $S^h$. For any vertex $z$ on $\mathcal{T}_h$, construct a
patch of elements $\mathcal{K}_{z}$ containing 
at least six vertices around $z$. Select all
vectices in $\mathcal{K}_z$ as sampling points and fit a quadratic polynomial
$p_z \in \mathbb{P}_2(\mathcal{K}_z)$ in least square sense at those
sampling points. Then the recovered gradient at $z$ is defined as
\begin{equation*}
    (G_hu_h)(z) = \nabla p_z(z).
\end{equation*}
$G_hu_h$  on the whole domain is obtained by interpolation. According to
\cite{zhangnaga2005, nagazhang2005}, $G_h$ enjoys the following properties
\begin{itemize}
    \item [(1)] $\|\nabla u - G_h u_I\|\lesssim h^2|u|_{3, \Omega}$,
	where $u_I$ is the linear interpolation of $u$ in $S_h$.
    \item [(2)] $\|G_hv_h\|_{0, \Omega}\lesssim \|\nabla v_h\| _{0, \Omega},
	\forall v_h \in S_h$.
\end{itemize}

According to \cite{nagazhangzhou2006}, two adjacent triangles (sharing
a common edge) form an $O(h^{1+\alpha})$ ($\alpha > 0$) approximate
parallelogram if the lengths of any two opposite edges differ by
only $O(h^{1+\alpha})$.

\begin{definition}
    The triangulation $\mathcal{T}_h$ is said to satisfy Condition
    $\alpha$ if any two adjacent triangles form an $O(h^{1+\alpha})$
    parallelogram.
\end{definition}

Using the same methods  \cite{zhangnaga2005, shen2006},
we can prove the following  superconvergence result:
\begin{theorem}
    Suppose $M(\lambda_i) \subset H_0^1(\Omega) \cap W^{3, \infty}(\Omega)$
    and  $\mathcal{T}_h$ satisfies Condition $\alpha$. Let $G_h$ be the
    polynomial preserving recovery operator. Then for any  eigenfunction
    of \eqref{equ:fem} corresponding to $\lambda_{i, h}$,
    there exists an eigenfunction $u_i\in M(\lambda_i)$
    corresponding to $\lambda_i$ such that
    \begin{equation}
	\|\mathcal{D}^{\frac{1}{2}}\nabla u_i - \mathcal{D}^{\frac{1}{2}}
	G_hu_{i, h}\|_{0, \Omega}\lesssim h^{1+\beta}
	\|u_i\|_{3, \infty, \Omega}, \quad \beta = \min(\alpha, 1).
	\label{equ:recovery}
    \end{equation}
    \label{thm:recovery}
\end{theorem}
As pointed out in \cite{nagazhangzhou2006}, $\alpha = \infty$ if
$\mathcal{T}_h$ is generated  using regular refinement. Fortunately,
the fine grid $\mathcal{T}_h$ is always a regular
refinement of some coarse grid $\mathcal{T}_H$ for two-grid method.
When we introduce two-grid methods in Section 3,
we only perform gradient recovery on fine grid $\mathcal{T}_h$ .
Thus we assume $\alpha = \infty$ and hence
$\beta = 1$ in section 3.

\section{Superconvergent two-grid methods} In the literature, two-grid methods
\cite{xuzhou2001, yangbi2011, hucheng}  were proposed to reduce the
 cost of eigenvalue computations. To further improve the accuracy,
two different approaches: gradient recovery enhancement
\cite{nagazhangzhou2006, shen2006, mengzhang2012} and two-space
methods \cite{andreev, racheva} can be used. Individually, those tools are
useful in certain circumstances.  Combined them properly, we are able to
design much effective and superconvergence algorithms, which we shall describe
below.

\subsection{Gradient recovery enhanced shifted inverse power two-grid scheme}
In this scheme, we first use the shifted inverse power based two-grid
scheme \cite{yangbi2011, hucheng} and then apply the gradient recovery
enhancing technique \cite{nagazhangzhou2006}.

\begin{algorithm}
\caption{}\label{alg:inverse}
\begin{enumerate}
\item  Solve  the eigenvalue problem on a coarse grid $\mathcal{T}_H$:
   Find $(\lambda_{i,H}, u_{i,H})\in \mathbb{R} \times S^H_0 $ and
 $\|u_{i, H}\|_a = 1$ satisfying
\begin{equation} \label{eqn:poisson}
  a(u_{i,H},v_H)=\lambda_{i,H} b(u_{i,H},v_H),\quad \forall v_H\in S^H_0.
\end{equation}
\item Solve a source problem on the fine grid $\mathcal{T}_h$:
    Find  $u^{i}_{h} \in S^h_0$
    such that
\begin{equation}
    a(u^{i}_{ h}, v_h)-\lambda_{i, H}(u^{i}_{h}, v_h) =
    (u_{i, H}, v_h), \quad \forall v_h\in S^{h}_0,
    \label{equ:inverse}
\end{equation}
and set $u^{i, h} = \frac{u^{i}_{ h}}{\|u^{i}_{ h}\|_{a}}$.
\item Apply the gradient recovery operator $G_h$ on $u^{i, h}$ to get
    $G_hu^{i, h}$.
\item Set
    \begin{equation}
	\lambda^{i, h} = \frac{a(u^{i, h}, u^{i, h})}{(u^{i, h}, u^{i, h})}
	- \frac{\|\mathcal{D}^{\frac{1}{2}}\nabla u^{i, h}-
	\mathcal{D}^{\frac{1}{2}}G_hu^{i, h}\|^2_{0, \Omega}}
	{(u^{i, h}, u^{i, h})}.
	\label{equ:improve}
    \end{equation}
\end{enumerate}
\end{algorithm}

In the proof of our main superconvergence result, we need
the following Lemma, which was proved in \cite[Theorem 4.1]{yangbi2011}.
\begin{lemma}
Suppose that $M(\lambda_i) \subset H^1_0(\Omega)\cap W^{3, \infty}(\Omega)$.
    Let $(\lambda^{i, h},
    u^{i, h})$ be an approximate eigenpair of \eqref{equ:variational} obtained
    by Algorithm \ref{alg:inverse} and let $H$ be properly small.
    Then
    \begin{equation}
	\text{dist}(u^{i, h}, M_{h}(\lambda_i))
	\lesssim H^4 + h^2,  \label{equ:inv}
    \end{equation}
    where $\text{dist}(u^{i, h}, M_{h}(\lambda_i))
    = \inf\limits_{v \in M_{h}(\lambda_i)}\|u^{i, h}- v\|_{a, \Omega}$.
    \label{lem:eigenfunction}
\end{lemma}

Based on the above Lemma, we can establish the superconvergence result for
eigenfunctions.
\begin{theorem}
    Suppose that $M(\lambda_i) \subset H^1_0(\Omega)\cap W^{3, \infty}(\Omega)$
    . Let $(\lambda^{i, h},
    u^{i, h})$ be an approximate eigenpair of \eqref{equ:variational} obtained
    by Algorithm \ref{alg:inverse} and let $H$ be properly small.
    Then there exists $u_i \in M(\lambda_i)$ such that
    \begin{equation}
	\|\mathcal{D}^{\frac{1}{2}}G_hu^{i, h} -
	\mathcal{D}^{\frac{1}{2}}\nabla u_i\|_{0, \Omega} \lesssim
	(H^4+h^2).
	\label{equ:inversefunerror}
    \end{equation}
    \label{thm:eigenfunction}
\end{theorem}
{\em Proof}. Let the
eigenfunctions $\{ u_{j, h}\}^{i+q-1}_{j=i}$  be an orthonormal basis of
$M_h(\lambda_i)$. Note that $$\text{dist}(u^{i, h}, M_{h}(\lambda_i))
=\|u^{i, h}-\sum_{j=i}^{j=i+q-1}a(u^{i, h},
u_{j, h})u_{j, h}\|_{a, \Omega}.$$
Let $\tilde{u}_h=\sum_{j=i}^{j=i+q-1}a(u^{i, h}, u_{j, h})u_{j, h}$.
 According to Theorem \ref{thm:recovery}, there
exist $\{\tilde{u}_j\}_{j=i}^{i+q-1} \subset M(\lambda_i)$ suth that
\begin{equation}
  \label{equ:temp}
    \|\mathcal{D}^{\frac{1}{2}}G_hu^{j, h} -
    \mathcal{D}^{\frac{1}{2}}\nabla \tilde{u}_j\|_{0, \Omega}
    \lesssim h^2.
  \end{equation}
Let $u_i = \sum_{j=i}^{j=i+q-1}a(u^{i, h}, u_{j, h})\tilde{u}_j$; then
$u_i\in M(\lambda_i)$. Using \eqref{equ:temp}, we can derive that
\begin{align*}
  &\|\mathcal{D}^{\frac{1}{2}}G_h\tilde{u}^{h} -
    \mathcal{D}^{\frac{1}{2}}\nabla {u}_i\|_{0, \Omega}\\
    =&\|\sum_{j=i}^{j=i+q-1}a(u^{i, h}, u_{j, h})
    (\mathcal{D}^{\frac{1}{2}}G_hu_{j, h} -
    \mathcal{D}^{\frac{1}{2}}\nabla \tilde{u}_j)\|_{0, \Omega}\\
    \lesssim &\left(\sum_{j=i}^{j=i+q-1}\|(\mathcal{D}
    ^{\frac{1}{2}}G_hu_{j, h} -
    \mathcal{D}^{\frac{1}{2}}\nabla \tilde{u}_j)\|^2_{0, \Omega}\right)
    ^{\frac{1}{2}}\\
    \lesssim &h^2.
  \end{align*}
  Thus, we have
  \begin{equation*}
  \begin{split}
    &\|\mathcal{D}^{\frac{1}{2}}G_hu^{i, h} -
    \mathcal{D}^{\frac{1}{2}}\nabla u_i\|_{0, \Omega}\\
    \le & \|\mathcal{D}^{\frac{1}{2}}G_h(u^{i, h} -
    \tilde{u}_h)\|_{0, \Omega} +
    \|\mathcal{D}^{\frac{1}{2}}G_h\tilde{u}_h-
   \mathcal{D}^{\frac{1}{2}}\nabla u_i\|_{0, \Omega}\\
   \lesssim & \|G_h(u^{i, h} -  \tilde{u}_h)\|_{0, \Omega} +h^2\\
   \lesssim & \|\nabla(u^{i, h} -  \tilde{u}_h)\|_{0, \Omega} + h^2\\
   \lesssim & \|u^{i, h} -  \tilde{u}_h\|_{a, \Omega} + h^2\\
   \lesssim & (H^4+h^2) + h^2\\
   \lesssim & H^4+h^2;
  \end{split}
  \end{equation*}
where we use  Lemma \ref{lem:eigenfunction} to bound
$\|u^{i, h} -  \tilde{u}_h\|_{a, \Omega}$. \endproof

The following Lemma is needed in the proof of  a  superconvergence property
of our eigenvalue approximation.
\begin{lemma}
    Suppose that $M(\lambda_i)\subset H_0^1(\Omega)\cap W^{3, \infty}(\Omega)$.
    Let $(\lambda^{i, h}, u^{i, h})$ be an approximate eigenpair
    of \eqref{equ:variational} obtained by
    Algorithm \ref{alg:inverse} and let H be properly small. Then
    \begin{equation}
	\|\mathcal{D}^{\frac{1}{2}}G_hu^{i, h} - \mathcal{D}^{\frac{1}{2}}
	\nabla  u^{i, h}\|_{0, \Omega}\lesssim (h+H^2).
	\label{equ:inverseeta}
    \end{equation}
    \label{lem:error}
\end{lemma}
{\em Proof}. Let $\tilde{u}_h$ be defined as in Theorem \ref{thm:eigenfunction}.
Then we have
\begin{equation*}
       \begin{split}
	 &\|\mathcal{D}^{\frac{1}{2}}G_hu^{i, h} -
	 \mathcal{D}^{\frac{1}{2}}\nabla u^{i, h}\|_{0, \Omega}  \\
	 \le & \|\mathcal{D}^{\frac{1}{2}}G_hu^{i, h} -
	 \mathcal{D}^{\frac{1}{2}}G_h \tilde{u}_h\|_{0, \Omega}
	 + \|\mathcal{D}^{\frac{1}{2}}G_h\tilde{u}_h -
	 \mathcal{D}^{\frac{1}{2}}\nabla \tilde{u}_h\|_{0, \Omega}
	 + \|\mathcal{D}^{\frac{1}{2}}\nabla \tilde{u}_h -
	 \mathcal{D}^{\frac{1}{2}}\nabla u^{i, h}\|_{0, \Omega}\\
       \lesssim & \|G_hu^{i, h} -
       G_h\tilde{u}_h\|_{0, \Omega}
	 + \|\mathcal{D}^{\frac{1}{2}}G_h\tilde{u}_h -
	 \mathcal{D}^{\frac{1}{2}}\nabla \tilde{u}_h\|_{0, \Omega}
	 + \|\mathcal{D}^{\frac{1}{2}}\nabla\tilde{ u}_h -
	 \mathcal{D}^{\frac{1}{2}}\nabla u^{i, h}\|_{0, \Omega}\\
	 \lesssim &\|\nabla u_{i, h} - \nabla \tilde{u}_h \| _{0, \Omega}
	 +\|\mathcal{D}^{\frac{1}{2}}G_h\tilde{u}_h -
	 \mathcal{D}^{\frac{1}{2}}\nabla \tilde{u}_h\|_{0, \Omega}\\
	 \lesssim &\|u_{i, h} - \tilde{u}_h \| _{a, \Omega}
	 +\|\mathcal{D}^{\frac{1}{2}}G_h\tilde{u}_h -
	 \mathcal{D}^{\frac{1}{2}}\nabla \tilde{u}_h\|_{0, \Omega}\\
	 \lesssim &(H^4+h^2) + h\\
	 \lesssim &(H^2+h).
       \end{split}
   \end{equation*}
   Here we use the fact that $\|\cdot\|_{a, \Omega}$ and $\|\cdot\|_{1, \Omega}$
   are two equivalent norms on $H^{1}_0(\Omega)$.\endproof

Now we are in a perfect position to prove our main  superconvergence result
for eigenvalue  approximation.
\begin{theorem}
     Suppose that $M(\lambda_i)\subset H_0^1(\Omega)\cap W^{3, \infty}(\Omega)$.
    Let $(\lambda^{i, h}, u^{i, h})$ be an approximate eigenpair
    of \eqref{equ:variational} obtained by
    Algorithm \ref{alg:inverse} and let H be properly small.
    \begin{equation}
	|\lambda^{i, h}-\lambda_i| \lesssim  H^6+h^3.
	\label{equ:inverseeigenvalueerror}
    \end{equation}
    \label{thm:eigenvalue}
\end{theorem}
{\em Proof}. It follows from \eqref{equ:identity} and \eqref{equ:improve} that

\begin{equation*}
    \begin{split}
	& \lambda^{i, h} - \lambda_i \\
	=&\frac{a(u^{i, h}, u^{i, h})}{(u^{i, h}, u^{i, h})}-
	\frac{\|\mathcal{D}^{\frac{1}{2}}\nabla u^{i, h} -
	\mathcal{D}^{\frac{1}{2}}G_h u^{i, h}\|^2_{0, \Omega}}
	{(u^{i, h}, u^{i, h})} -\lambda_i \\
	=&\frac{a(u^{i, h}-u_i, u^{i, h} - u_i)}{(u^{i, h}, u^{i, h})}
	-\frac{\|\mathcal{D}^{\frac{1}{2}}\nabla u^{i, h} -\mathcal{D}^{\frac{1}{2}}G_h u^{i, h}\|^2_{0, \Omega}}
	{(u^{i, h}, u^{i, h})}	
	-\frac{\lambda_i(u^{i, h} - u_i, u^{i, h} - u_i)}
	{(u^{i, h}, u^{i, h})}\\
	=&\frac{(\mathcal{D}^{\frac{1}{2}}(u^{i, h}-u_i),
	\mathcal{D}^{\frac{1}{2}}(u^{i, h} - u_i))}{(u^{i, h}, u^{i, h})}
	-\frac{\|\mathcal{D}^{\frac{1}{2}}\nabla u^{i, h} -
	\mathcal{D}^{\frac{1}{2}}G_h u^{i, h}\|^2_{0, \Omega}}
	{(u^{i, h}, u^{i, h})}	+\\
	&\frac{ (c(u^{i, h} -u_i), u^{i, h} - u_i) -
	\lambda_i(u^{i, h} - u_i, u^{i, h} - u_i)}
	{(u^{i, h}, u^{i, h})}\\
	= &\frac{\|\mathcal{D}^{\frac{1}{2}}\nabla G_h u^{i, h} -
	\mathcal{D}^{\frac{1}{2}}\nabla u_i\|^2_{0, \Omega}}
	{(u^{i, h}, u^{i, h})}
	+\frac{2(\mathcal{D}^{\frac{1}{2}}G_h u^{i, h} -
	\mathcal{D}^{\frac{1}{2}}\nabla u_i,
	\mathcal{D}^{\frac{1}{2}}\nabla u^{i, h} -
	\mathcal{D}^{\frac{1}{2}} G_h u^{i, h})}
	{(u^{i, h}, u^{i, h})}+\\
	&\frac{ (c(u^{i, h} -u_i), u^{i, h} - u_i) -
	\lambda_i(u^{i, h} - u_i, u^{i, h} - u_i)}
	{(u^{i, h}, u^{i, h})}.
    \end{split}
\end{equation*}
From Theorem $4.1$ in \cite{yangbi2011}, we know that
$\|u^{i, h} - u_i\|_{0, \Omega}\lesssim (H^4+h^2)$ and hence
the last term in the above equation is bounded by $O((H^4+h^2)^2)$.
Theorem \ref{thm:eigenfunction} implies that the first term is also
bounded by $O((H^4+h^2)^2)$. Using the H\"{o}lder inequality,  we obtain
\begin{eqnarray}\label{holder}
    && |(\mathcal{D}^{\frac{1}{2}}G_h u^{i, h} - \mathcal{D}^{\frac{1}{2}}\nabla u_i,
    \mathcal{D}^{\frac{1}{2}}\nabla u^{i, h} -
    \mathcal{D}^{\frac{1}{2}}G_h u^{i, h})| \nonumber\\
    &\le&\|\mathcal{D}^{\frac{1}{2}}G_h u^{i, h} -
    \mathcal{D}^{\frac{1}{2}}\nabla u_i\|_{0, \Omega}
    \|\mathcal{D}^{\frac{1}{2}}\nabla u^{i, h} -
    \mathcal{D}^{\frac{1}{2}}G_h u^{i, h})\|_{0, \Omega} \nonumber\\
    &\lesssim &(H^4+h^2)(H^2+h) \lesssim  H^6+h^3
\end{eqnarray}
and hence
\begin{equation*}
    |\lambda^{i, h} - \lambda_i | \lesssim  H^6+h^3.
\end{equation*}
This completes our proof.\endproof

Taking $H = O(\sqrt{h})$, Theorem \ref{thm:eigenfunction} and
\ref{thm:eigenvalue} implies that we can get $O(h^2)$ superconvergence
and $O(h^3)$ superconvergence for eigenfunction and eigenvalue
approximation,
respectively.

\begin{remark}
\label{rmk:better}
Using the H\"{o}lder inequality to estimate (\ref{holder})
does not take into account the cancellation in the integral.
Similar as \cite{nagazhangzhou2006}, numerical experiments show that the
 actual  bound is
\begin{equation*}
    |(\mathcal{D}^{\frac{1}{2}}G_h u^{i, h} -
    \mathcal{D}^{\frac{1}{2}}\nabla u_i,
    \mathcal{D}^{\frac{1}{2}}\nabla u^{i, h} -
    \mathcal{D}^{\frac{1}{2}}G_h u^{i, h})|
    \lesssim (H^4+h^2)^2,
\end{equation*}
\end{remark}
which says that we have ``double''-order gain by applying recovery.


\begin{remark}
    Algorithm \ref{alg:inverse} is a combination of
    the shifted inverse power two-grid method \cite{yangbi2011, hucheng}
     and gradient recovery enhancement \cite{nagazhangzhou2006}. It inherits all excellent
    properties of both methods: low computational cost and superconvergence.
     We will demonstrate in our numerical tests that Algorithm 1 outperforms
    shifted inverse power two-grid method in \cite{yangbi2011, hucheng}.
\end{remark}

\begin{remark}
  If we firstly use classical two-grid methods as in \cite{xuzhou2001} and then
  apply gradient recovery, we can  prove
  $\|\mathcal{D}^{\frac{1}{2}}G_hu^{i, h} -
	\mathcal{D}^{\frac{1}{2}}\nabla u_i\|_{0, \Omega} \lesssim
	(H^2+h^2)$  and
    $|\lambda^{i, h} - \lambda_i | \lesssim  H^3+h^3$.
    It means we can only get optimal convergence rate insteading of 
    superconvergent convergence rate when $H=O(\sqrt{h})$.
\end{remark}

\subsection{Higher order space based superconvergent two-grid scheme}
Our second scheme can be viewed as a combination of the two-grid scheme
proposed by Yang and Bi \cite{yangbi2011} or Hu and Cheng \cite{hucheng}
and the two-space method introduced by Racheva and Andreev \cite{racheva}.

\begin{algorithm}
\caption{}\label{alg:high}
\begin{enumerate}
\item  Solve an eigenvalue problem on a coarse grid $\mathcal{T}_H$:
   Find $(\lambda_{i,H}, u_{i,H})\in \mathbb{R} \times S^H_0 $ and
 $\|u_{i, H}\|_a = 1$ satisfying
\begin{equation} \label{eqn:sourse}
  a(u_{i,H},v_H)=\lambda_{i,H} b(u_{i,H},v_H),\quad \forall v_H\in S^H_0.
\end{equation}

\item Solve a source problem on the fine grid $\mathcal{T}_h$:
    Find  $u^{i}_{h} \in S^{h, 2}_0$
    such that
\begin{equation}
    a(u^{i,h}, v_h)-\lambda_{i, H}(u^{i, h}, v_h) =
    (u_{i, H}, v_h), \quad \forall v_h\in S^{h, 2}_0.
    \label{equ:high}
\end{equation}
\item Compute the Rayleigh quotient
    \begin{equation}
	\lambda^{i, h} = \frac{a(u^{i, h}, u^{i, h})}{(u^{i, h}, u^{i, h})}.
	\label{equ:ray}
    \end{equation}
\end{enumerate}
\end{algorithm}

Note that we use linear finite element space $S^H_0$ on coarse grid
$\mathcal{T}_H$ and quadratic finite element space $S^{h, 2}_0$ on
fine grid $\mathcal{T}_h$. Compared with the two-grid scheme
\cite{yangbi2011, hucheng},
the main difference is that Algorithm \ref{alg:high} uses linear
element on coarse grid $\mathcal{T}_H$ and quadratic element on
fine grid $\mathcal{T}_h$ while the two-grid uses linear element
on both coarse grid $\mathcal{T}_H$ and  $\mathcal{T}_h$.
Compared with the two-space method \cite{racheva}, the main difference
is that Algorithm 2 uses a coarse grid $\mathcal{T}_H$
and a fine grid $\mathcal{T}_h$ whereas the two-space method only uses a
grid $\mathcal{T}_h$.  Algorithm \ref{alg:high} shares the advantages
of both methods: low computational cost and high accuracy. Thus, we would
expect Algorithm \ref{alg:high} performs much better than both methods.

For Algorithm 2, we have the following Theorem:
\begin{theorem}
    Suppose that $M(\lambda_i)\subset H_0^1(\Omega)\cap H^{3}(\Omega)$.
    Let $(\lambda^{i, h}, u^{i, h})$ be an approximate eigenpair
    of \eqref{equ:variational} by
    Algorithm \ref{alg:inverse} and let H be properly small.
    Then there exists $u_i \in M(\lambda_i)$ such that
    \begin{equation}
      |u^{i, h}-u_i|_{a, \Omega} \lesssim (H^4+h^2);
	    \label{equ:highfunction}
	\end{equation}
    \begin{equation}
         \lambda^{i, h} - \lambda_i \lesssim (H^8+h^4).
	\label{equ:highvalue}
    \end{equation}
    \label{thm:high}
\end{theorem}
{\em Proof}. By Theorem $4.1$ in \cite{yangbi2011}, we have
 \begin{equation}
     |u^{i, h}-u_i|_{a, \Omega} \lesssim \eta_{a}(H)\delta_H^3(\lambda_i)
     + \delta_{h}(\lambda_i);
	    \label{equ:yangfunction}
	\end{equation}
and
    \begin{equation}
         \lambda^{i, h} - \lambda_i \lesssim \eta_{a}^2(H)\delta_H^6(\lambda_i)
     + \delta_{h}^2(\lambda_i).
	\label{equ:yangvalue}
    \end{equation}
Since we use linear element on $\mathcal{T}_H$ and quadratic element on
$\mathcal{T}_h$, it follows from the interpolation error estimate
\cite{brenner, ciarlet} that
\begin{equation*}
    \eta_a(H)\lesssim H, \quad \delta_H(\lambda_i)\lesssim H,\quad
    \delta_h(\lambda_i) \lesssim h^2.
\end{equation*}
Substituting the above three estimate into \eqref{equ:yangfunction} and
\eqref{equ:yangvalue}, we get \eqref{equ:highfunction} and
\eqref{equ:highvalue}.\endproof

Comparing Algorithm \ref{alg:inverse} and \ref{alg:high}, the main difference
 is that Algorithm \ref{alg:inverse} solves a source problem on
fine grid $\mathcal{T}_h$ using linear element and hence perform gradient
recovery while Algorithm \ref{alg:high} solves a source problem on fine
grid $\mathcal{T}_h$ using quadratic element. Both Algorithm \ref{alg:inverse}
and \ref{alg:high} lead to $O(h^2)$ superconvergence for eigenfunction
approximation and $O(h^4)$ ultraconvergence for eigenvalue approximation by
taking $H=O(\sqrt{h})$. The message we would like to deliver here is that
polynomial preserving recovery plays a similar role as quadratic element,
but with much lower computational cost.

\begin{remark}
    In order to get higher order convergence, we require higher regualrity
    such as $M(\lambda_i) \subset H_0^1(\Omega)\cap W^{3, \infty}(\Omega)$ for
    Algorithm \ref{alg:inverse} and
     $M(\lambda_i) \subset H_0^1(\Omega)\cap H^{3}(\Omega)$ for
     Algorithm \ref{alg:high}, in the proof.
     However, we can use Algorithm \ref{alg:inverse}
     and \ref{alg:high} to get high accuracy approximation even
     with low regularity.
\end{remark}

\section{Multilevel adaptive methods}

\begin{algorithm}
\caption{Given a tolerance $\epsilon >0$ and a parameter $0\le \theta < 1$.}
\label{alg:adaptive}
\begin{enumerate}
    \item Generate an initial mesh $\mathcal{T}_{h_0}$.
    \item Solve \eqref{equ:variational} on $\mathcal{T}_{h_{0}}$ to get a discrete
	eigenpair
	$(\bar{\lambda}^{h_0}, u^{h_0})$.
    \item Set $\ell =0$.
    \item Compute $\eta(u^{h_{\ell}}, T)$ and
      $\eta(u^{h_{\ell}}, \Omega)$, then let
	\begin{equation*}
	    \lambda^{h_{\ell}} =  \bar{\lambda}^{h_{\ell}}-
	    \eta(u^{h_{\ell}}, \Omega)^2.
	\end{equation*}
    \item If $\eta(u^{h_{\ell}}, \Omega)^2<\epsilon$, stop; else go to 6.
    \item Choose a minimal subset  of elements
	$\widehat{\mathcal{T}}_{h_{\ell}}\subset \mathcal{T}_{h_{\ell}}$
	such that
	\begin{equation*}
	    \sum_{T\in \widehat{\mathcal{T}}_{h_{\ell}}}
	         \eta^2(u_h, T)\ge \theta
	     \eta^2(u_h, \Omega);
	\end{equation*}
	then refine the elements in $\widehat{\mathcal{T}}_{h_{\ell}} $ and
	 necessary elements
	to get a new conforming mesh $\mathcal{T}_{h_{\ell+1}}$.
    \item Find $u\in S^{h_{\ell+1}}_0$ such that
	\begin{equation*}
	    a(u, v) = \lambda_{h_{\ell}} b(u^{h_{\ell}}, v),
	    \quad v \in S^{h_{\ell+1}}_0,
	\end{equation*}
	and set $u^{h_{\ell+1}} = \frac{u}{\|u\|_{0, \Omega}}$. Define
	\begin{equation}
	    \bar{\lambda}^{h_{\ell+1}}
	    = \frac{a(u^{h_{\ell+1}}, u^{h_{\ell+1}})}
	    {b(u^{h_{\ell+1}}, u^{h_{\ell+1}})}.
	    \label{equ:bvp1}
	\end{equation}
    \item Let $\ell = \ell +1$ and go to 4.
\end{enumerate}
\end{algorithm}

\begin{algorithm}[!h]
\caption{Given a tolerance $\epsilon >0$ and a parameter $0\le \theta < 1$.}
\label{alg:improve}
\begin{enumerate}
    \item Generate an initial mesh $\mathcal{T}_{h_0}$.
    \item Solve \eqref{equ:variational} on $\mathcal{T}_{h_{0}}$ to get a discrete
	eigenpair
	$(\bar{\lambda}^{h_0}, u^{h_0})$.
    \item Set $\ell =0$.
  \item Compute $\eta(u^{h_{\ell}}, T)$ and
      $\eta(u^{h_{\ell}}, \Omega)$, then let
	\begin{equation*}
	    \lambda^{h_{\ell}} =  \bar{\lambda}^{h_{\ell}}-
	    \eta(u^{h_{\ell}}, \Omega)^2.
	\end{equation*}
    \item If $\eta(u^{h_{\ell}}, \Omega)^2<\epsilon$, stop; else go to 6.
\item Choose a minimal subset  of elements
	$\widehat{\mathcal{T}}_{h_{\ell}}\subset \mathcal{T}_{h_{\ell}}$
	such that
	\begin{equation*}
	    \sum_{T\in \widehat{\mathcal{T}}_{h_{\ell}}}
	         \eta^2(u_h, T)\ge \theta
	     \eta^2(u_h, \Omega);
	\end{equation*}
	then refine the elements in $\widehat{\mathcal{T}}_{h_{\ell}} $ and
	 necessary elements
	to get a new conforming mesh $\mathcal{T}_{h_{\ell+1}}$.
    \item Find $u\in S^{h_{\ell+1}}_0$ such that
	\begin{equation}
	    a(u, v) - \lambda_{h_{\ell}} b(u, v) = b(u^{h_{\ell}}, v)
	    , \quad v \in S^{h_{\ell+1}}_0,
	    \label{equ:bvp2}
	\end{equation}
	and set $u^{h_{\ell+1}} = \frac{u}{\|u\|_{0, \Omega}}$. Define
	\begin{equation*}
	    \bar{\lambda}^{h_{\ell+1}}
	    = \frac{a(u^{h_{\ell+1}}, u^{h_{\ell+1}})}
	    {b(u^{h_{\ell+1}}, u^{h_{\ell+1}})}.
	\end{equation*}
    \item Let $\ell = \ell +1$ and go to 4.
\end{enumerate}
\end{algorithm}

In this section, we incorporate two-grid methods and gradient recovery
enhancing technique into the framework of adaptive finite element method
and propose two multilevel adaptive methods.  Both methods
only need to solve an eigenvalue problem on initial mesh and solve
an associated boundary value problem on adaptive refined mesh
during every iteration.

Let $u_h$ be a finite element solution in $S^{h}$ and $G_h$ be PPR
recovery operator. Define a local a posteriori error estimator on the element
$T$ as
\begin{equation}
    \eta(u_h, T) = \|\mathcal{D}^{\frac{1}{2}}G_hu_h -
    \mathcal{D}^{\frac{1}{2}}\nabla u_h\|_{0, T},
    \label{equ:indicator}
\end{equation}
and a global error estimator as
\begin{equation}
    \eta(u_h, \Omega) = \left(\sum_{T\in \mathcal{T}_h}
    \eta(u_h, T)\right)^{\frac{1}{2}}.
    \label{equ:globalindicator}
\end{equation}
Given a tolerance $\epsilon$ and a parameter $\theta$, we describe our
multilevel adaptive methods in Algorithm \ref{alg:adaptive} and
\ref{alg:improve}. Here we use D$\ddot{\text{o}}$rfler marking strategy
\cite{dorfler} in step 6.

Note that the only difference between  Algorithm \ref{alg:adaptive}
and \ref{alg:improve} is that they solve different boundary value problems
on step 7. Algorithm \ref{alg:adaptive} solves  boundary value problem
\eqref{equ:bvp1} like two-grid scheme in \cite{xuzhou2001} while
Algorithm \ref{alg:improve} solves boundary value problem \eqref{equ:bvp2}
similar to two-grid scheme in \cite{yangbi2011, hucheng}.
Boundary value problem \eqref{equ:bvp2} would lead to a near singular linear
system. Although there are many efficient iterative methods, like multigrid
methods, as pointed out in \cite{hucheng}, the computational cost of
solving \eqref{equ:bvp1} should be higher than solving \eqref{equ:bvp2}.
Numerical results of both methods are almost the same as indicated by
examples in next section. Thus, Algorithm \ref{alg:adaptive} is highly
recommended.

Compared to methods in \cite{liyang, xie}, Algorithm \ref{alg:adaptive}
and \ref{alg:improve} use recovery based a  posteriori error estimator.
The propose of gradient recovery in the above two algorithms is twofold.
The first one is  to provide an  asymptotically exact a posteriori error
estimator.
The other is to greatly improve  the accuracy of eigenvalue and eigenfunction
approximations.
Superconvergence result $O(N^{-1})$ and ultraconvergence $O(N^{-2})$ are
numerically observed for eigenfunction and eigenvalue approximation
respectively. However, methods in \cite{liyang, xie} can only numerically give
asymptotically optimal results.
We want to emphasize that the new algorithms  can get
superconvergence or ultraconvergence
results with no more or even less computational cost compared to
the methods proposed in \cite{liyang, xie}.

\section{Numerical Experiment}
In this section, we present several numerical examples to demonstrate
the effectiveness  and  superconconvergence of the proposed algorithms
and validity  our theoretical results. All algorithms are implemented
using finite element package iFEM developed by Chen \cite{chen}.

The first  example is designed to demonstrate superconvergence
property of Algorithm \ref{alg:inverse} and  \ref{alg:high} and make
some comparison with the two-grid scheme in \cite{yangbi2011, hucheng}.
Let the $i$th
eigenpairs obtained by Algorithm \ref{alg:inverse} and
\ref{alg:high} be denoted by 
$(\lambda^{i,\text{A1}}, u^{i,\text{A1}})$ and
$(\lambda^{i,\text{A2}}, u^{i,\text{A2}})$. Also, let
$(\lambda^{i,\text{TG}}, u^{i,\text{TG}})$ be the  $i$th eigenpair produced by
the shift inverse based two-grid scheme in \cite{yangbi2011, hucheng}.

The presentation of other examples are to illustrate the effectiveness and
sueprconvergence of Algorithm \ref{alg:adaptive} and  \ref{alg:improve}.
In these examples, we focus on the first eigenpair. 
Let $\bar{\lambda}_{\text{A3}}$ and $\lambda_{\text{A3}}$ be
the eigenvalue generated by Algorithm \ref{alg:adaptive} without and
with gradient recovery enhancing, respectively. Define
$\bar{\lambda}_{\text{A4}}$, $\lambda_{\text{A4}}$, $u_{\text{A3}}$,
and $u_{\text{A4}}$ in a similar way.

{\bf Example 1.}  Consider the following Laplace eigenvalue problem
\begin{equation}
    \begin{cases}
	-\Delta u = \lambda u, \quad \text{in } \Omega, \\
	u = 0 , \quad \text{on } \partial \Omega,
    \end{cases}
    \label{equ:laplacesquare}
\end{equation}
where $\Omega= (0, 1)\times (0, 1)$. The eigenvalue of \eqref{equ:laplacesquare}
are $\lambda_{k, l} = (k^2+l^2)\pi^2$ and  the corresponding eigenfunctions
are $u_{k, l} = \sin(k\pi)\sin(l\pi)$ with $k, l = 1, 2, \cdots$.
It is easy to see the first three eigenvalues are $\lambda_1 = 2\pi^2$
and $\lambda_2=\lambda_3=5\pi^2$.

\begin{figure}[!h]
  \centering
  \begin{minipage}[c]{0.5\textwidth}
  \centering
  \includegraphics[width=0.9\textwidth]{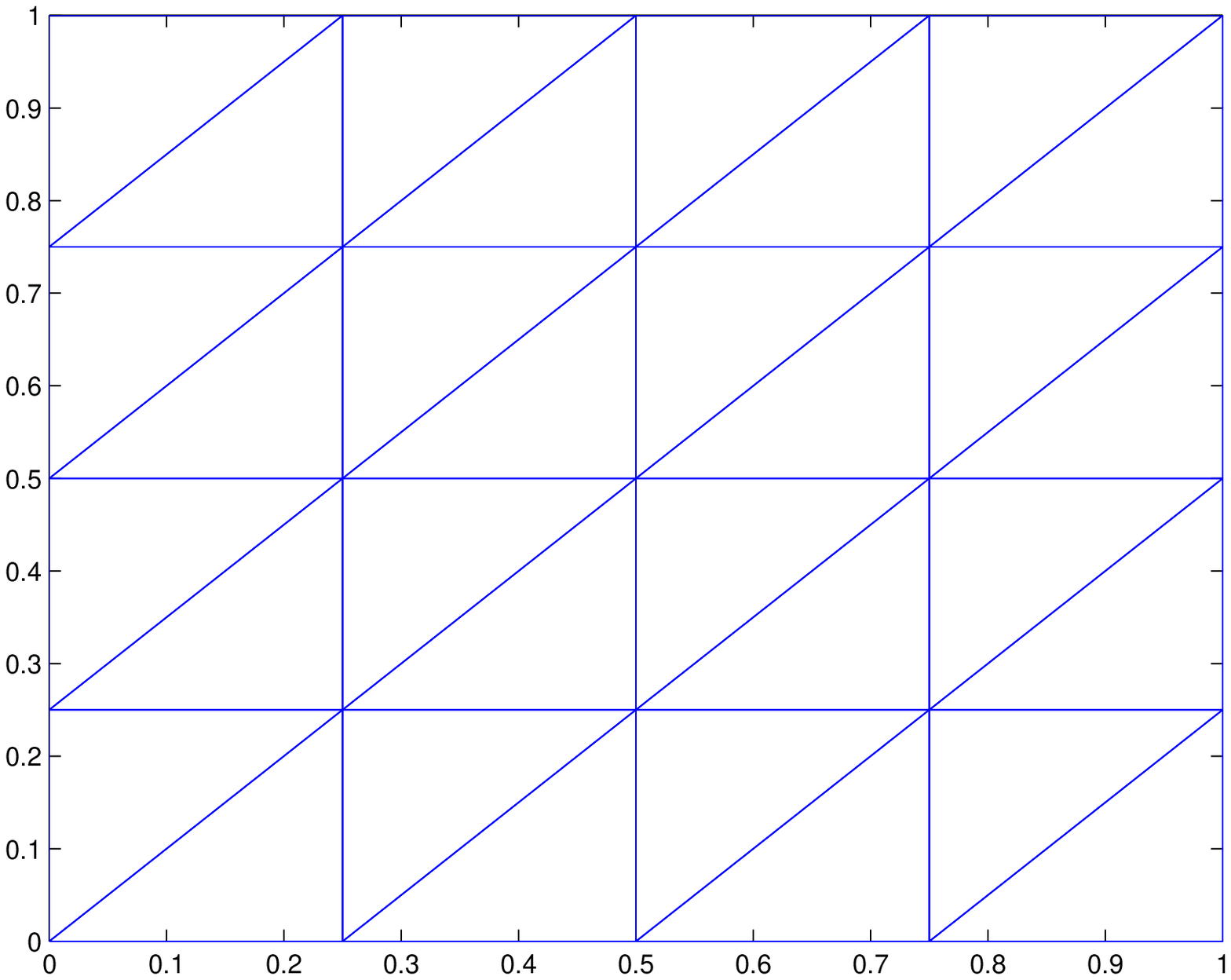}
  \caption{niform Mesh for Example 1}
\label{fig:ex1uniform}
\end{minipage}%
 \begin{minipage}[c]{0.5\textwidth}
  \centering
  \includegraphics[width=0.9\textwidth]{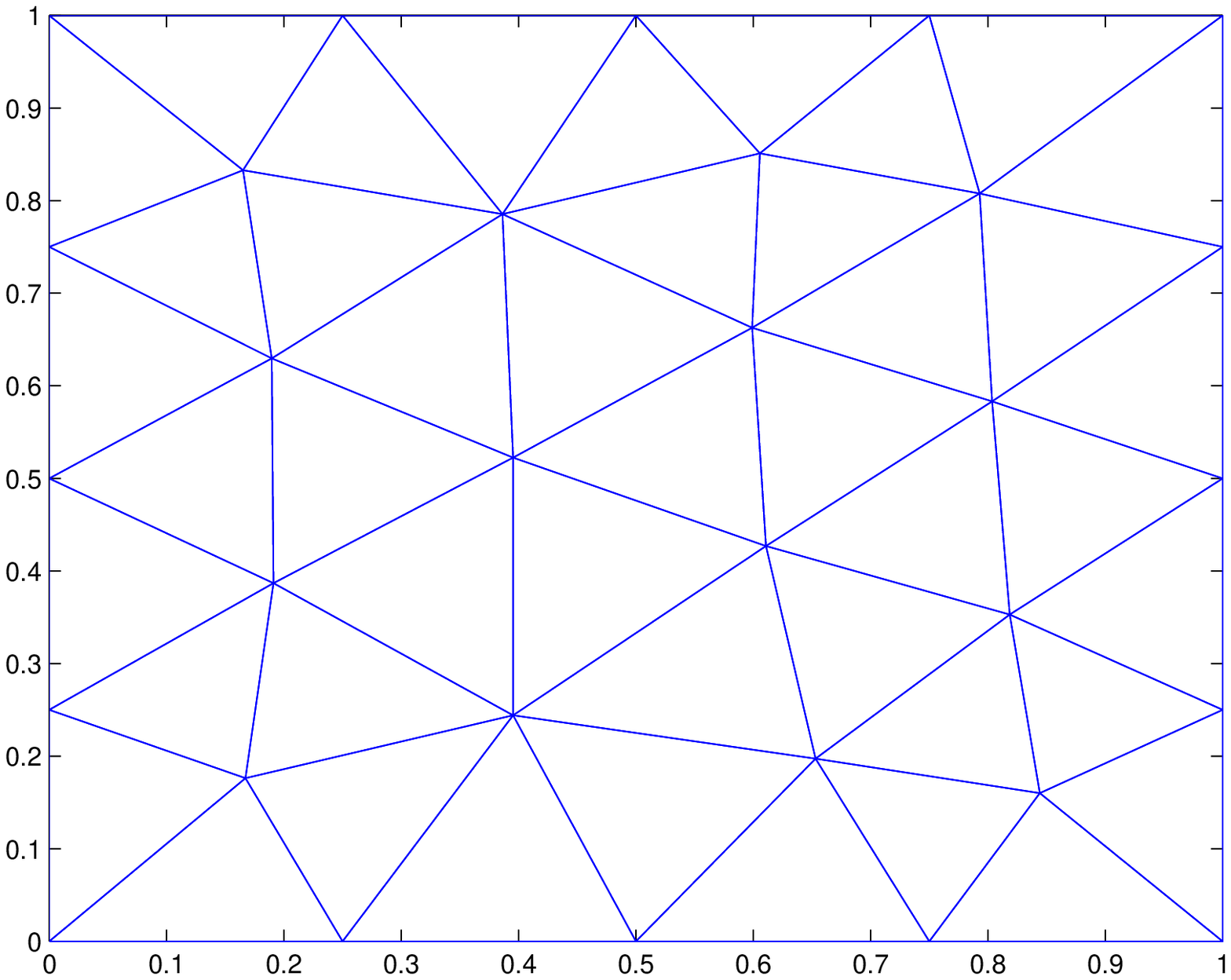}
  \caption{Delaunay Mesh for Example 1}
\label{fig:ex1delaunay}
\end{minipage}
\end{figure}

\begin{table}[htb!]
\centering
\footnotesize
\caption{ Eigenpair errors of 
  Algorithm \ref{alg:inverse}  for Example 1 on Uniform Mesh}
\begin{tabular}{|c|c|c|c|c|c|c|c|}
\hline
i&H&h&$\lambda^{i,\text{A1}}$ & $\lambda^{i,\text{A1}}-\lambda_i$ &
$Order$ & $\|G_hu^{i,\text{A1}} - \nabla u_i\|_{0, \Omega}$&$Order$\\ \hline
1&1/4&1/16&19.733813512912&-5.40e-03&&7.059395e-02&\\ \hline
1&1/8&1/64&19.739186935311&-2.19e-05&3.97&4.387700e-03&2.00\\ \hline
1&1/16&1/256&19.739208716241&-8.59e-08&4.00&2.734342e-04&2.00\\ \hline
1&1/32&1/1024&19.739208801843&-3.36e-10&4.00&1.707544e-05&2.00\\ \hline
2&1/4&1/16&49.311524605286&-3.65e-02&0.00&------------&\\ \hline
2&1/8&1/64&49.347897768530&-1.24e-04&4.10&------------&\\ \hline
2&1/16&1/256&49.348021565420&-4.40e-07&4.07&------------&\\ \hline
2&1/32&1/1024&49.348022003783&-1.66e-09&4.02&------------&\\ \hline
3&1/4&1/16&49.311750580349&-3.63e-02&0.00&------------&\\ \hline
3&1/8&1/64&49.347802761238&-2.19e-04&3.69&------------&\\ \hline
3&1/16&1/256&49.348021182216&-8.23e-07&4.03&------------&\\ \hline
3&1/32&1/1024&49.348022002296&-3.15e-09&4.01&------------&\\ \hline
\end{tabular}\label{tab:ex1uniforma1}
\end{table}

\begin{table}[htb!]
\centering
\footnotesize
\caption{ Eigenpair errors of 
  Algorithm \ref{alg:high}  for Example 1 on Uniform Mesh}
\begin{tabular}{|c|c|c|c|c|c|c|c|}
\hline
i&H & h & $\lambda^{i,\text{A2}}$ & $\lambda^{i,\text{A2}}-\lambda_i$ &
$Order$ & $\|\nabla u^{i,\text{A2}} - \nabla u_i\|_{0, \Omega}$&$Order$\\ \hline
1&1/4&1/16&19.740140941323&9.32e-04&&3.344371e-02&\\ \hline
1&1/8&1/64&19.739212357340&3.56e-06&4.02&2.076378e-03&2.00\\ \hline
1&1/16&1/256&19.739208816236&1.41e-08&3.99&1.308168e-04&1.99\\ \hline
1&1/32&1/1024&19.739208802235&5.59e-11&3.99&8.198527e-06&2.00\\ \hline
2&1/4&1/16&49.399143348018&5.11e-02&0.00&------------&\\ \hline
2&1/8&1/64&49.348217238157&1.95e-04&4.02&------------&\\ \hline
2&1/16&1/256&49.348022827362&8.22e-07&3.95&------------&\\ \hline
2&1/32&1/1024&49.348022008741&3.29e-09&3.98&------------&\\ \hline
3&1/4&1/16&49.573605264596&2.26e-01&0.00&------------&\\ \hline
3&1/8&1/64&49.348559514553&5.38e-04&4.36&------------&\\ \hline
3&1/16&1/256&49.348024046492&2.04e-06&4.02&------------&\\ \hline
3&1/32&1/1024&49.348022013418&7.97e-09&4.00&------------&\\\hline
\end{tabular}
\label{tab:ex1uniforma2}
\end{table}

\begin{table}[htb!]
\centering
\footnotesize
\caption{ Eigenpair errors of 
 shift-inverse  Two-grid scheme
 for Example 1 on Uniform Mesh}
\begin{tabular}{|c|c|c|c|c|c|c|c|}
\hline
i& H & h & $\lambda^{i,\text{TG}}$ & $\lambda^{i,\text{TG}}-\lambda_i$ &
$Order$ & $\|\nabla u^{i,\text{TG}} - \nabla u_i\|_{0, \Omega}$&$Order$\\ \hline
1&1/4&1/16&19.930259632276&1.91e-01&&4.375101e-01&\\ \hline
1&1/8&1/64&19.751103117985&1.19e-02&2.00&1.090672e-01&1.00\\ \hline
1&1/16&1/256&19.739951989101&7.43e-04&2.00&2.726155e-02&1.00\\ \hline
1&1/32&1/1024&19.739255250511&4.64e-05&2.00&6.815303e-03&1.00\\ \hline
2&1/4&1/16&50.199210624678&8.51e-01&0.00&-----------&\\ \hline
2&1/8&1/64&49.399315353599&5.13e-02&2.03&-----------&\\ \hline
2&1/16&1/256&49.351217793553&3.20e-03&2.00&-----------&\\ \hline
2&1/32&1/1024&49.348221696982&2.00e-04&2.00&-----------&\\ \hline
3&1/4&1/16&50.779973345337&1.43e+00&0.00&-----------&\\ \hline
3&1/8&1/64&49.428220994371&8.02e-02&2.08&-----------&\\ \hline
3&1/16&1/256&49.353003975409&4.98e-03&2.00&-----------&\\ \hline
3&1/32&1/1024&49.348333256327&3.11e-04&2.00&-----------&\\ \hline
\end{tabular}
\label{tab:ex1uniformtg}
\end{table}

First, uniform mesh as in  Fig \ref{fig:ex1uniform} is considered.
The fine meshes $\mathcal{T}_h$ are of sizes $h = 2^{-j}$ ($j = 4, 6, 8, 10$)
and the corresponding  coarse meshes $\mathcal{T}_H$ of size $H=\sqrt{h}$.
 Table \ref{tab:ex1uniforma1} lists the numerical results for
Algorithm \ref{alg:inverse}. 
$\|G_h u^{i,\text{A1}}-\nabla u_i\|_{0, \Omega}$
(i = 1)
superconverges at rate of $O(h^2)$ which consists with our theoretical
analysis. However, $|\lambda^{\text{i,A1}}-\lambda_i|$ 
(i = 1, 2, 3)
ultraconverges at rate of $O(h^4)$
which is better than the results predicted by Theorem \ref{thm:eigenvalue}.
In particular, it verifies the statement  in Remark
\ref{rmk:better}.
Since $\lambda_2$ and $\lambda_2$ are multiples eigenvalues, the error 
of eigenfunctions approximation are not available and it is 
represented by $-$ in Tables \ref{tab:ex1uniforma1}-\ref{tab:ex1delaunaytg}.
One important thing we want to point out is that
we observe numerically that $\lambda_{\text{A1}}$ obtained by
Algorithm \ref{alg:inverse} approximates the exact eigenvalue from below; see
column $4$ in Table \ref{tab:ex1uniforma1}. Similar phenomenon was observed
in \cite{fang2013}  where they use a local high-order interpolation recovery.
We want to remark that lower bound of eigenvalue is very important in
practice and there are many efforts are made to obtain eigenvalue approximation
from below. The readers are referred to
\cite{duran2004, lin2014, yang2010, zhang2007} for other ways to approxiate
eigenvalue from below. In Table \ref{tab:ex1uniforma2}, we report the numerical
result of Algorithm \ref{alg:high}. As expected, $O(h^4)$ convergence of
eigenvalue approximation and $O(h^2)$ convergence of eigenfunction
approximation are observed which validate our Theorem \ref{thm:high}.
The   shift-inverse power method based 
two-grid scheme in \cite{yangbi2011, hucheng} is then considered,
the result being displayed in Table \ref{tab:ex1uniformtg}.
$\lambda^{i,\text{TG}}$ approximates $\lambda_i$ (i = 1, 2, 3)  at a
rate $O(h^2)$ and $\|u^{i,\text{TG}}-u_i\|_{a, \Omega}$ (i=1) converges at a
rate of $O(h)$.

Comparing Tables \ref{tab:ex1uniforma1} to \ref{tab:ex1uniformtg},
huge advantages of Algorithm \ref{alg:inverse} and \ref{alg:high}
are demonstrated. For instance, on the fine grid with size $h = 1/1024$
and corresponding coarse grid with size $H= 1/32$,
the  approximate  first eigenvalues
produced by Algorithm \ref{alg:inverse} and \ref{alg:high} are exact
up to $10$ digits while one can only trust the first five digits
of the first eigenvalue generated by the two-grid scheme in
\cite{yangbi2011, hucheng}.

\begin{table}[htb!]
\centering
\footnotesize
\caption{Comparison of Three Algorithms for 
Example 1 on Uniform mesh}
\begin{tabular}{|c|c|c|c|c|c|c|c|c|}
\hline
i&H & h & $i,\lambda^{\text{A1}}$& $\lambda^{\text{i,A1}}-\lambda_i$&
$\lambda^{i,\text{A2}}$&  $\lambda^{i,\text{A2}}-\lambda_i$&
$\lambda^{i,\text{TG}}$& $\lambda^{i,\text{TG}}-\lambda_i$ \\ \hline
1&1/2 & 1/16 & 20.1083669 & 3.69e-01& 20.2080796 & 4.69e-01 & 20.3504780 & 6.11e-01  \\ \hline
1&1/4 & 1/256 & 19.7398503 & 6.41e-04& 19.7398588 & 6.50e-04 & 19.7406011 & 1.39e-03  \\ \hline
\end{tabular}
\label{tab:compare}
\end{table}

Then we consider the case $H=O(\sqrt[4]{h})$  for the first 
eigenvalue. We use the fine meshes
of mesh size $h = 2^{-j}$ with $j = 4, 8$ and corresponding coarse
meshes satisfying $H=\sqrt[4]{h}$. The numerical results are showed in
Table \ref{tab:compare}. We can see that the two proposed  Algorithms
give   better approximate eigenvalues.  Thus Algorithm \ref{alg:inverse}
and \ref{alg:high} outperforms the two-grid scheme even in the case
$H=\sqrt[4]{h}$.  
One interesting thing that we want to mention is that
$\lambda^{i,\text{A1}}$ approximates $\lambda_i$ from above in this case,
see column $4$ in Table \ref{tab:compare}.

\begin{table}[htb!]
\centering
\footnotesize
\caption{ Eigenpair errors of 
  Algorithm \ref{alg:inverse}  for Example 1 on Delaunay Mesh}
\begin{tabular}{|c|c|c|c|c|c|c|c|}
\hline
i&H&h&$\lambda^{i,\text{A1}}$ & $\lambda^{i,\text{A1}}-\lambda_i$ &
$Order$ & $\|G_hu^{i,\text{A1}} - \nabla u_i\|_{0, \Omega}$&$Order$\\ \hline
1&31&385&19.735647110619&-3.56e-03&&5.338236e-02&\\ \hline
1&105&5761&19.739198229599&-1.06e-05&2.15&2.835582e-03&1.08\\ \hline
1&385&90625&19.739208765246&-3.69e-08&2.05&1.686396e-04&1.02\\ \hline
1&1473&1443841&19.739208802041&-1.38e-10&2.02&1.049196e-05&1.00\\ \hline
2&31&385&49.307472112236&-4.05e-02&0.00&------------&\\ \hline
2&105&5761&49.347888708818&-1.33e-04&2.11&------------&\\ \hline
2&385&90625&49.348021524994&-4.80e-07&2.04&------------&\\ \hline
2&1473&1443841&49.348022003630&-1.82e-09&2.01&------------&\\ \hline
3&31&385&49.301142920140&-4.69e-02&0.00&------------&\\ \hline
3&105&5761&49.347856273486&-1.66e-04&2.09&------------&\\ \hline
3&385&90625&49.348021393237&-6.12e-07&2.03&------------&\\ \hline
3&1473&1443841&49.348022003123&-2.32e-09&2.01&------------&\\ \hline
\end{tabular}
\label{tab:ex1delaunaya1}
\end{table}

\begin{table}[htb!]
\centering
\footnotesize
\caption{ Eigenpair errors of 
  Algorithm \ref{alg:high}  for Example 1 on Delaunay Mesh}
\begin{tabular}{|c|c|c|c|c|c|c|c|c|}
\hline
i&H & h & $\lambda^{i,\text{A2}}$ & $\lambda^{i,\text{A2}}-\lambda_i$ &
$Order$ & $\|\nabla u^{i,\text{A2}} - \nabla u_i\|_{0, \Omega}$&$Order$\\ \hline
1&31&385&19.739293668773&8.49e-05&&9.258930e-03&\\ \hline
1&105&5761&19.739209125443&3.23e-07&2.06&5.705799e-04&1.03\\ \hline
1&385&90625&19.739208803434&1.26e-09&2.01&3.555028e-05&1.01\\ \hline
1&1473&1443841&19.739208802184&5.33e-12&1.97&2.220103e-06&1.00\\ \hline
2&31&385&49.350648806465&2.63e-03&0.00&------------&\\ \hline
2&105&5761&49.348029138391&7.13e-06&2.18&------------&\\ \hline
2&385&90625&49.348022031328&2.59e-08&2.04&------------&\\ \hline
2&1473&1443841&49.348022005547&1.00e-10&2.01&------------&\\ \hline
3&31&385&49.351570779092&3.55e-03&0.00&------------&\\ \hline
3&105&5761&49.348029733509&7.73e-06&2.27&------------&\\ \hline
3&385&90625&49.348022033250&2.78e-08&2.04&------------&\\ \hline
3&1473&1443841&49.348022005554&1.07e-10&2.01&------------&\\ \hline
\end{tabular}
\label{tab:ex1delaunaya2}
\end{table}

\begin{table}[htb!]
\centering
\footnotesize
\caption{ Eigenpair errors of 
 shift-inverse  Two-grid scheme
 for Example 1 on Delaunay Mesh}
\begin{tabular}{|c|c|c|c|c|c|c|c|c|}
\hline
i& H & h & $\lambda^{i,\text{TG}}$ & $\lambda^{i,\text{TG}}-\lambda_i$ &
$Order$ & $\|\nabla u^{i,\text{TG}} - \nabla u_i\|_{0, \Omega}$&$Order$\\ \hline
1&31&385&19.821235920927&8.20e-02&&2.865766e-01&\\ \hline
1&105&5761&19.744334806708&5.13e-03&1.02&7.159881e-02&0.51\\ \hline
1&385&90625&19.739529185236&3.20e-04&1.01&1.789929e-02&0.50\\ \hline
1&1473&1443841&19.739228826191&2.00e-05&1.00&4.474820e-03&0.50\\ \hline
2&31&385&49.828430094852&4.80e-01&0.00&------------&\\ \hline
2&105&5761&49.377951127988&2.99e-02&1.03&------------&\\ \hline
2&385&90625&49.349892261888&1.87e-03&1.01&------------&\\ \hline
2&1473&1443841&49.348138895061&1.17e-04&1.00&------------&\\ \hline
3&31&385&49.893495693695&5.45e-01&0.00&------------&\\ \hline
3&105&5761&49.381970792689&3.39e-02&1.03&------------&\\ \hline
3&385&90625&49.350143791388&2.12e-03&1.01&------------&\\ \hline
3&1473&1443841&49.348154618353&1.33e-04&1.00&------------&\\ \hline
\end{tabular}\label{tab:ex1delaunaytg}
\end{table}

Now, we turn to unstructured meshes. First we generate a coarse mesh
$\mathcal{T}_H$ and repeat regular refinement on $\mathcal{T}_H$
until $H=O(\sqrt{h})$ to get the corresponding fine mesh $\mathcal{T}_h$.
The first level coarse mesh is generated by EasyMesh \cite{easymesh}
and the other three level coarse mesh are generated by regular refinement.
The numerical results are provided in Tables \ref{tab:ex1delaunaya1} to
\ref{tab:ex1delaunaytg}. Note that $N_H$ and $N_h$ denote the number of
vertices on coarse mesh $\mathcal{T}_H$ and fine mesh $\mathcal{T}_h$,
respectively. Concerning the convergence of eigenvalue,
Algorithm \ref{alg:inverse} and \ref{alg:high}  ultraconverge at
rate $O(h^4)$ while the two-grid scheme converges at rate $O(h^2)$.
Note that in Tables 5.5--5.7, $N_H \approx H^{-2}$ and $N_h \approx h^{-2}$.
Therefore, convergent rates for $H$ and $h$ ``double" the rates for $N_H$ and $N_h$, respectively.
As for eigenfunction,  $\|G_h u^{i,\text{A1}}-\nabla u_i\|_{0, \Omega}$
and $\|\nabla u^{i,\text{A2}}-\nabla u_i\|_{0, \Omega}$  are about
$O(h^2)$ while $\|\nabla u^{\text{i,TG}}-\nabla u_i\|_{0, \Omega}\approx O(h)$.

{\bf Example 2.} In the previous example, the eigenfunctions $u$ are analytic.
Here we consider Laplace eigenvalue value problem  on
the L-shaped domain $\Omega = (-1, 1)\times (-1, 1)/ [0, 1) \times (-1, 0]$.
The first eigenfunction has a singularity at the origin. To capture  this
singularity, multilevel adaptive algorithms \ref{alg:adaptive} and
\ref{alg:improve} are used with $\theta = 0.4$.
Since the first exact eigenvalue is not available, we choose  an
approximation  $\lambda = 9.6397238440219$
obtained by Betcke and Trefethen in \cite{betcke},
which is correct up to $14$ digits.

\begin{figure}[!h]
  \centering
  \begin{minipage}[c]{0.5\textwidth}
  \centering
  \includegraphics[width=0.9\textwidth]{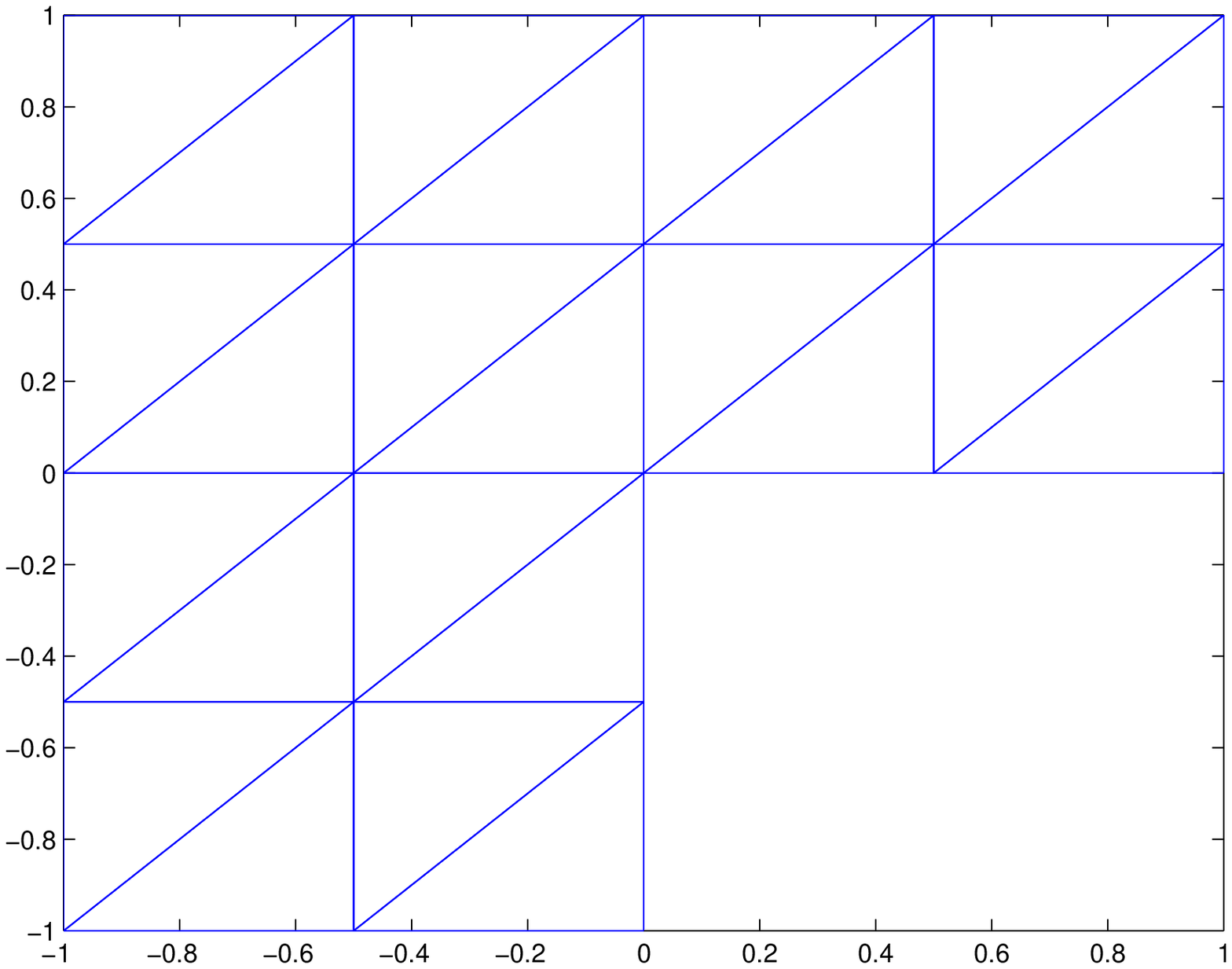}
  \caption{Initial Mesh for Example 2}
\label{fig:Lshape_initial}
\end{minipage}%
 \begin{minipage}[c]{0.5\textwidth}
  \centering
  \includegraphics[width=0.9\textwidth]{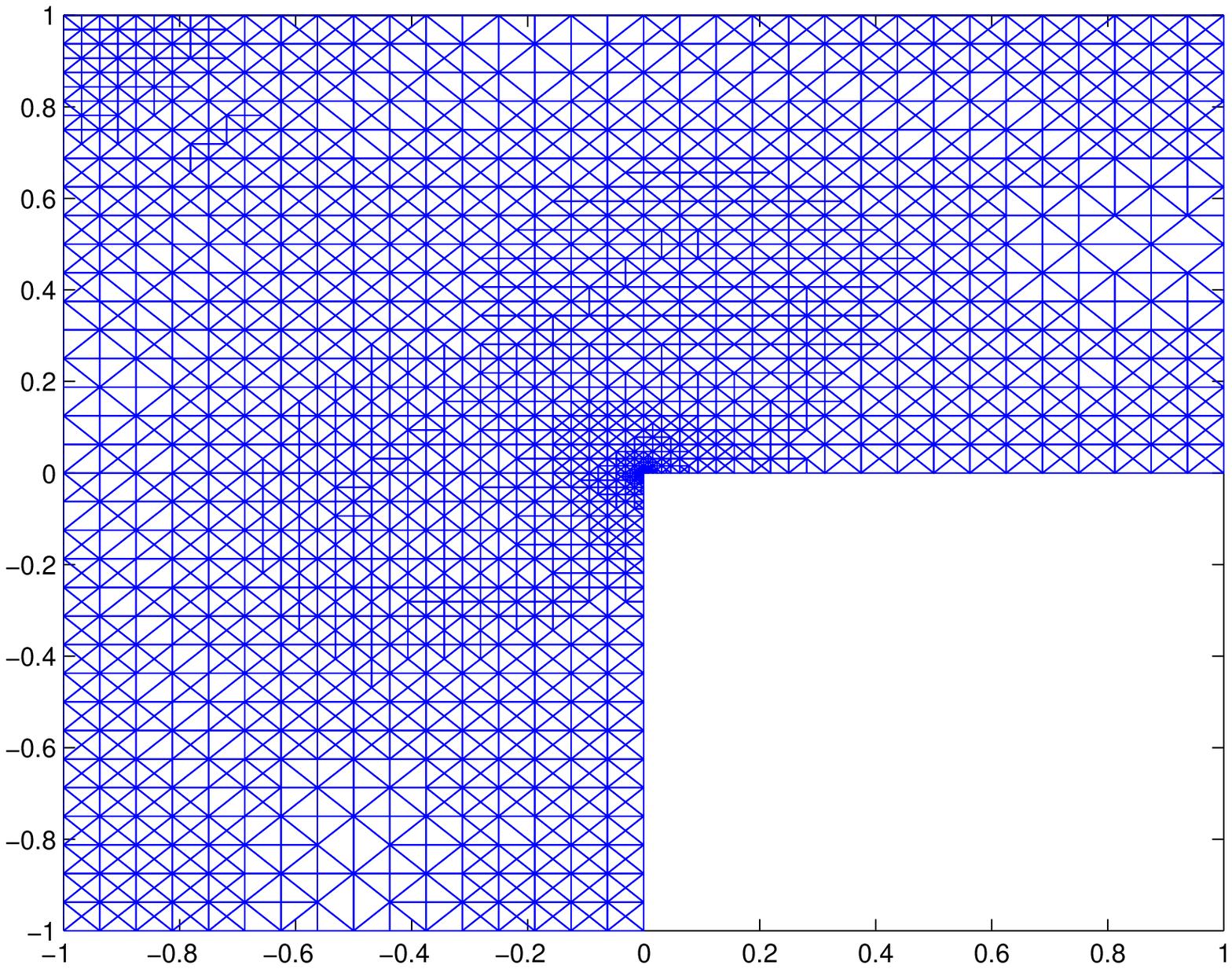}
  \caption{Adaptive Mesh for Example 2}
\label{fig:Lshape_adaptive}
\end{minipage}
\end{figure}

\begin{figure}[!h]
  \centering
  \includegraphics[width=0.8\textwidth]{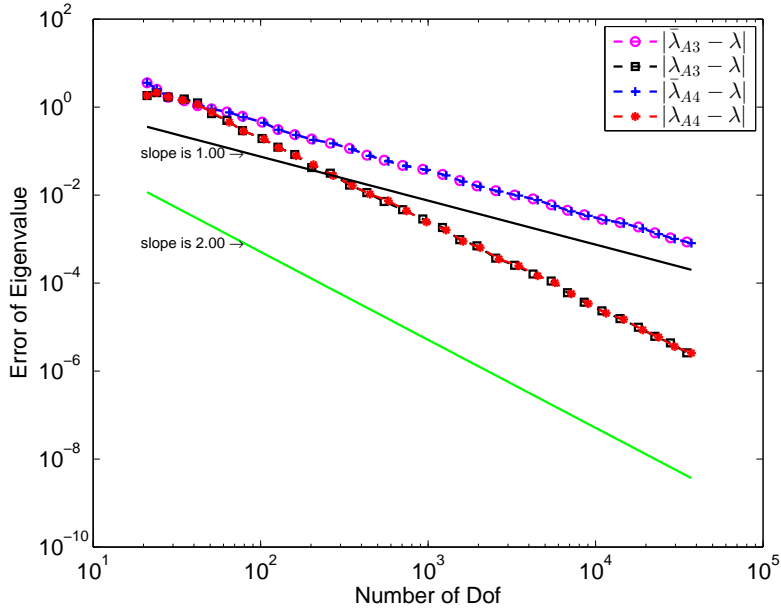}
  \caption{Eigenvalue Approximation Error  for Example 2}
\label{fig:Lshape_error}
\end{figure}

Fig \ref{fig:Lshape_initial}  shows the initial uniform mesh while
Fig \ref{fig:Lshape_adaptive} is the mesh after 18 adaptive iterations.
Fig \ref{fig:Lshape_error} reports numerical results of the first
eigenvalue approximation. It indicates clearly $\bar{\lambda}_{\text{A3}}$
and $\bar{\lambda}_{\text{A4}}$ approximate $\lambda$ at a rate of $O(N^{-1})$
while $\lambda_{\text{A3}}$ and $\lambda_{\text{A4}}$ approximate
$\lambda$ at a rate of $O(N^{-2})$.  The numerical results for
Algorithms \ref{alg:adaptive} and \ref{alg:improve} are almost the same.
Furthermore,  we notice that $\lambda_{\text{A3}}$
and $\lambda_{\text{A4}}$ approximate the exact eigenvalue from below.
It is well known that $\bar{\lambda}_{\text{A3}}$ and $\bar{\lambda}_{\text{A4}}$
are upper bounds for the exact eigenvalue. In actual computation,
we use $\bar{\lambda}_{\text{A3}}-\lambda_{\text{A3}}\le \epsilon$
as stop criteria for adaptive Algorithm  \ref{alg:adaptive} where
$\epsilon$ is the given tolerance.  A similar procedure is applied to
Algorithm \ref{alg:improve}.

\begin{figure}[!h]
  \centering
  \begin{minipage}[c]{0.5\textwidth}
  \centering
  \includegraphics[width=0.9\textwidth]{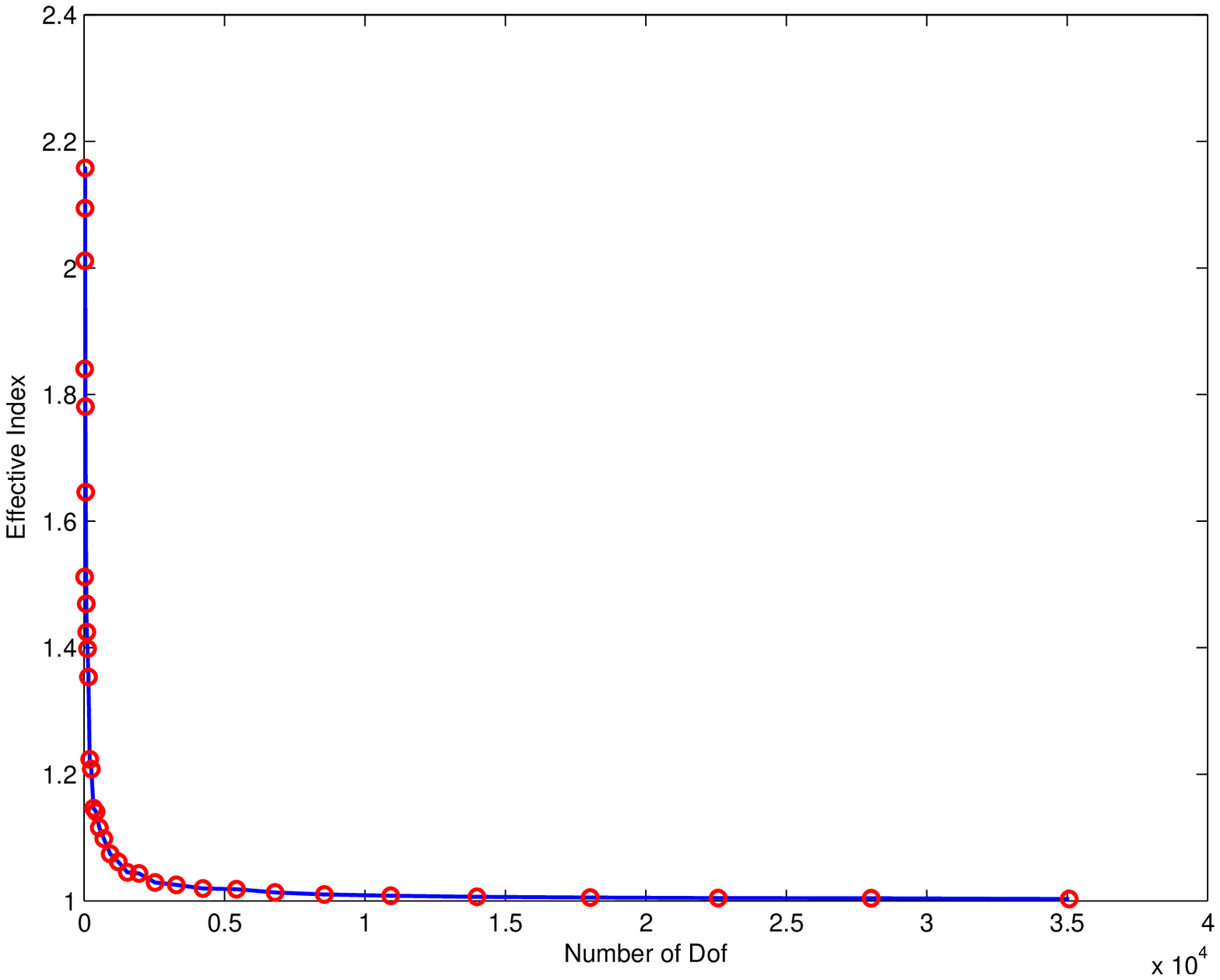}
  \caption{Effective index of Algorithm \ref{alg:adaptive} for Example 2}
\label{fig:effectindexa3lshape}
\end{minipage}%
 \begin{minipage}[c]{0.5\textwidth}
  \centering
  \includegraphics[width=0.9\textwidth]{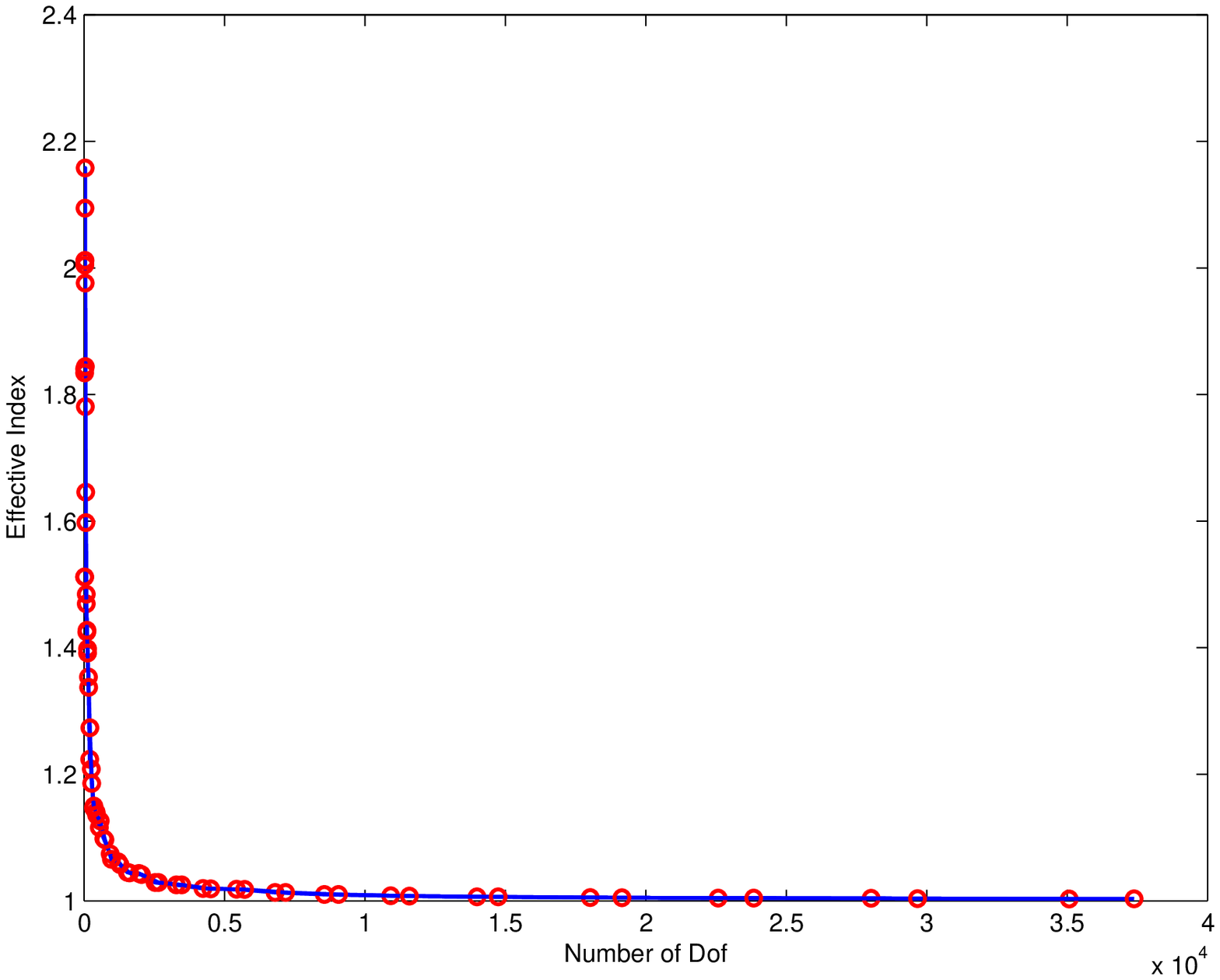}
  \caption{Effective index of Algorithm \ref{alg:improve} for Example 2}
\label{fig:effectindexa4lshape}
\end{minipage}
\end{figure}

In the context of adaptive finite element method for boundary value problems,
the effectivity index $\kappa$ is used to measure the quality of
an error estimator \cite{oden, babuska}. 
For eigenvalue problem,
it is better to consider eigenvalue effectivity  index
insteading of traditional effectivity index in \cite{oden, babuska}.
In the article, we consider  a similiar eigenvalue effective index
as in \cite{giani}
\begin{equation}
    \kappa = \frac{\|\mathcal{D}^{\frac{1}{2}}
    G_h u_h - \mathcal{D}^{\frac{1}{2}} \nabla u_h\|^2_{0, \Omega}}
    {|\lambda - \lambda_h|},
    \label{equ:effect}
\end{equation}
where $u_h$ is either $u_{\text{A3}}$ or $u_{\text{A4}}$ 
and $\lambda_h$ is either $\lambda_{\text{A3}}$ or $\lambda_{\text{A4}}$.
The effectivity index for the two proposed multilevel adaptive algorithms
are reported in  Figs \ref{fig:effectindexa3lshape} and
\ref{fig:effectindexa4lshape}.
 We see that $\kappa$ converges to $1$ quickly after the first few iterations, which
indicates that the posteriori error estimator \eqref{equ:indicator} or
\eqref{equ:globalindicator} is asymptotically exact.

{\bf Example 3.} Consider the following harmonic oscillator equation
\cite{greiner}, which
is a simple model in quantum mechanics,
\begin{equation}
      -\frac{1}{2}\Delta u + \frac{1}{2}|x|^2u = \lambda u,
      \quad \text{in } \mathbb{R}^2,
    \label{equ:electric}
\end{equation}
where $|x|=\sqrt{|x_1|^2+|x_2|^2}$. The first eigenvalue of
\eqref{equ:electric} is $\lambda = 1$ and the corresponding eigenfunction
is $u = \gamma e^{-|x|^2/2}$ with any nonzero constant $\gamma$.

\begin{figure}[!h]
  \centering
  \begin{minipage}[c]{0.5\textwidth}
  \centering
  \includegraphics[width=0.9\textwidth]{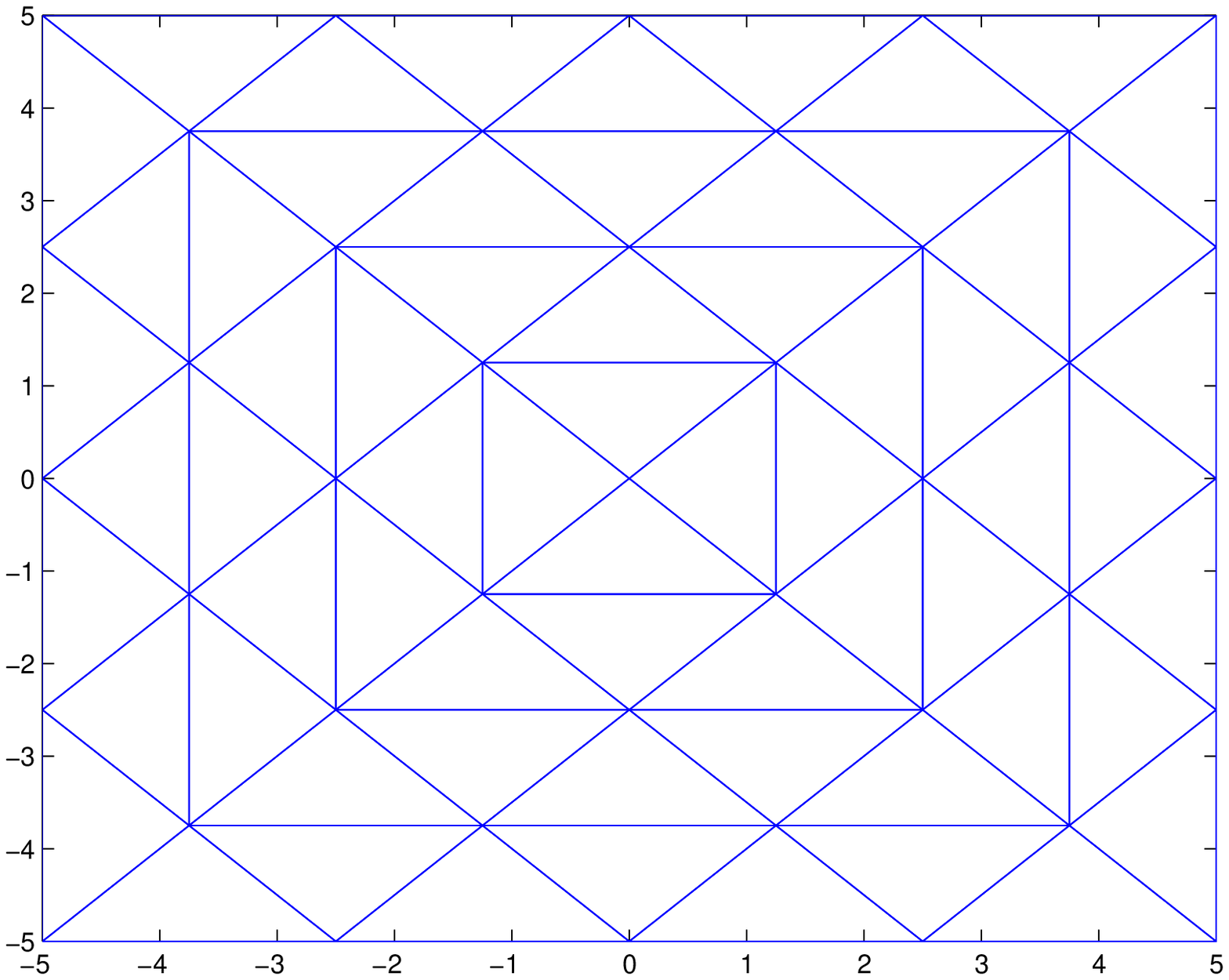}
  \caption{Initial Mesh for Example 3}
\label{fig:Electric_initial}
\end{minipage}%
 \begin{minipage}[c]{0.5\textwidth}
  \centering
  \includegraphics[width=0.9\textwidth]{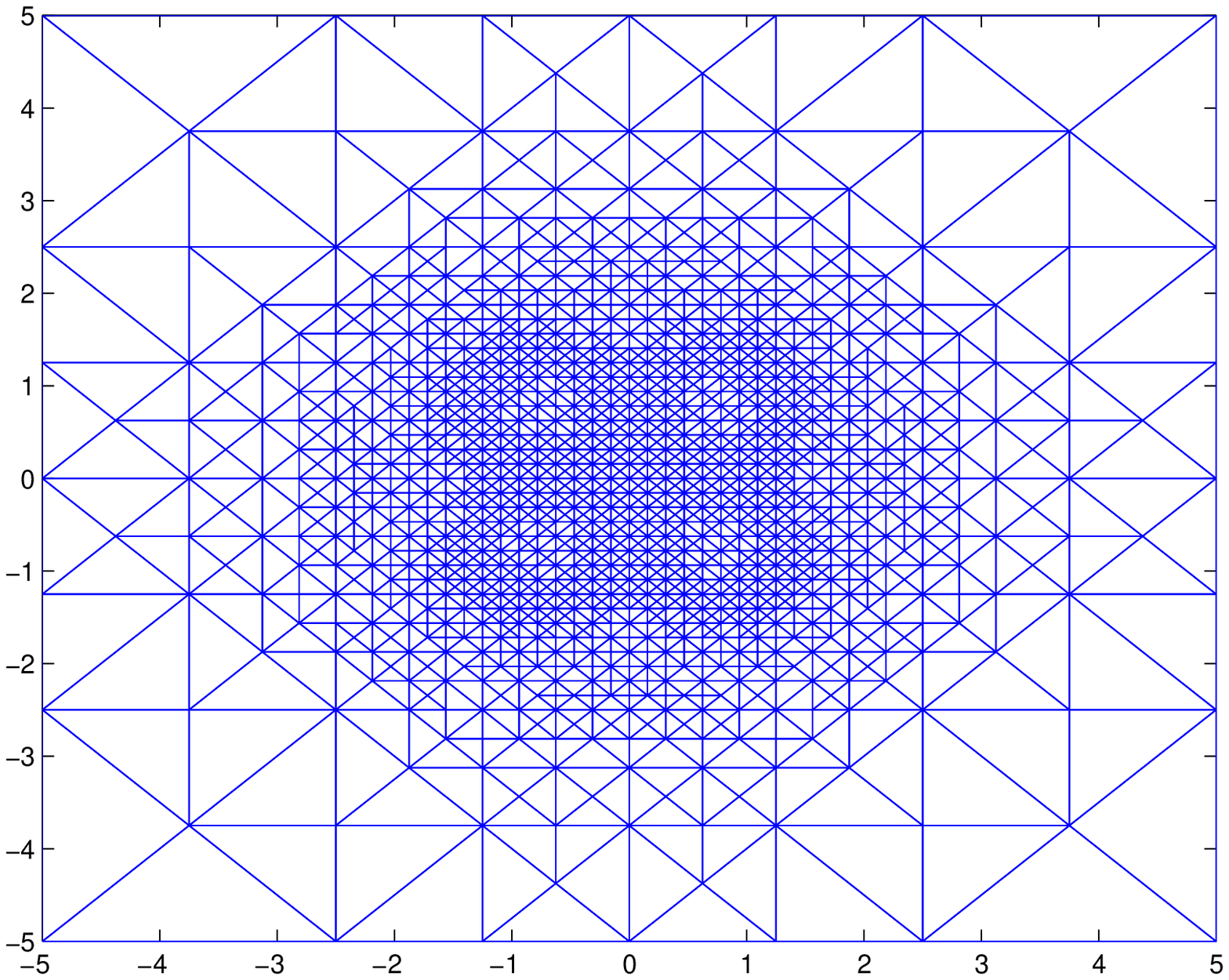}
  \caption{Adaptive Mesh for Example 3}
  \label{fig:Electric_adaptive}
\end{minipage}
\end{figure}

\begin{figure}[!h]
  \centering
  \begin{minipage}[c]{0.5\textwidth}
  \centering
  \includegraphics[width=0.9\textwidth]{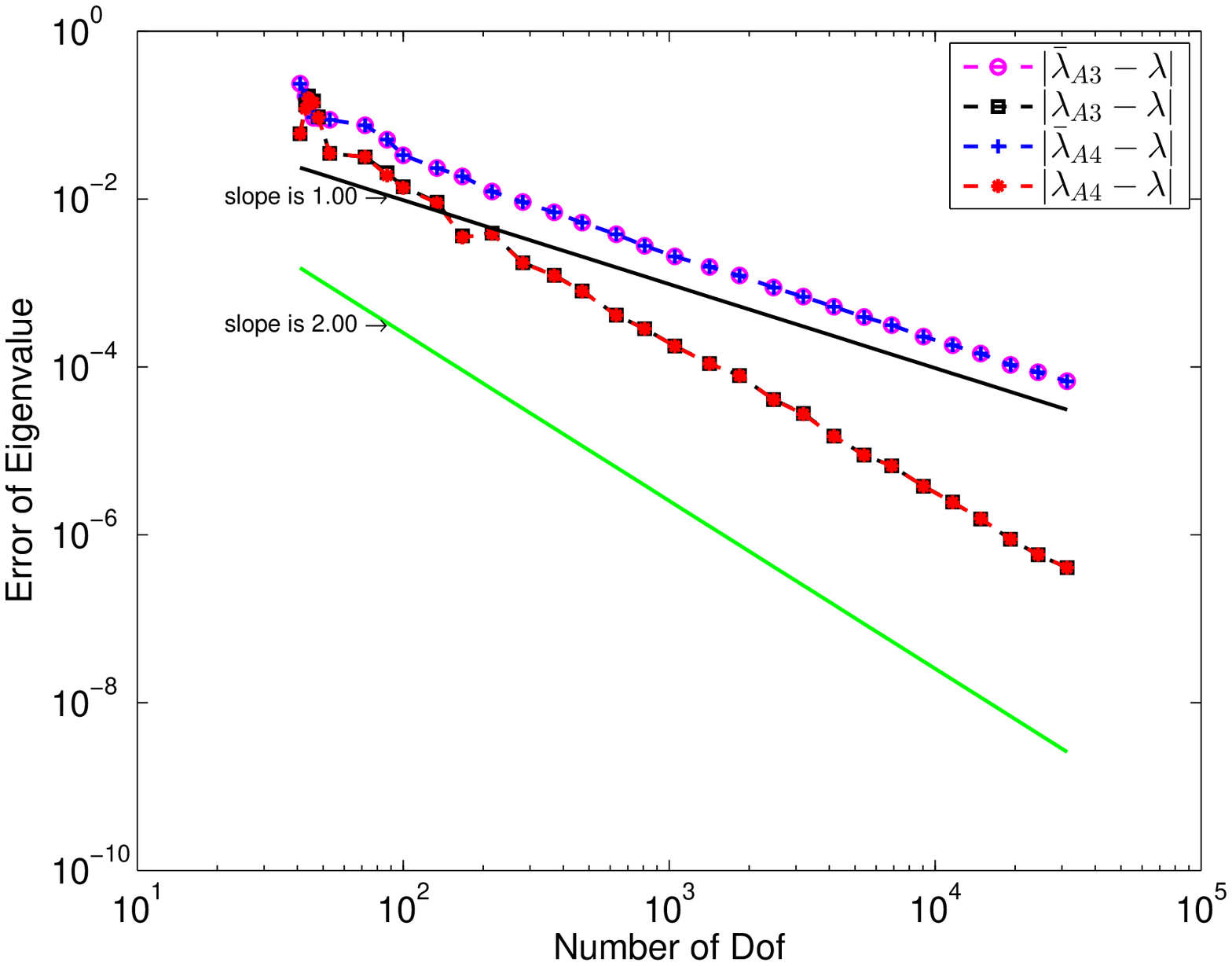}
  \caption{Eigenvalue approximatio Error for Example 3}
\label{fig:Electric_eigenvalue}
\end{minipage}%
 \begin{minipage}[c]{0.5\textwidth}
  \centering
  \includegraphics[width=0.9\textwidth]{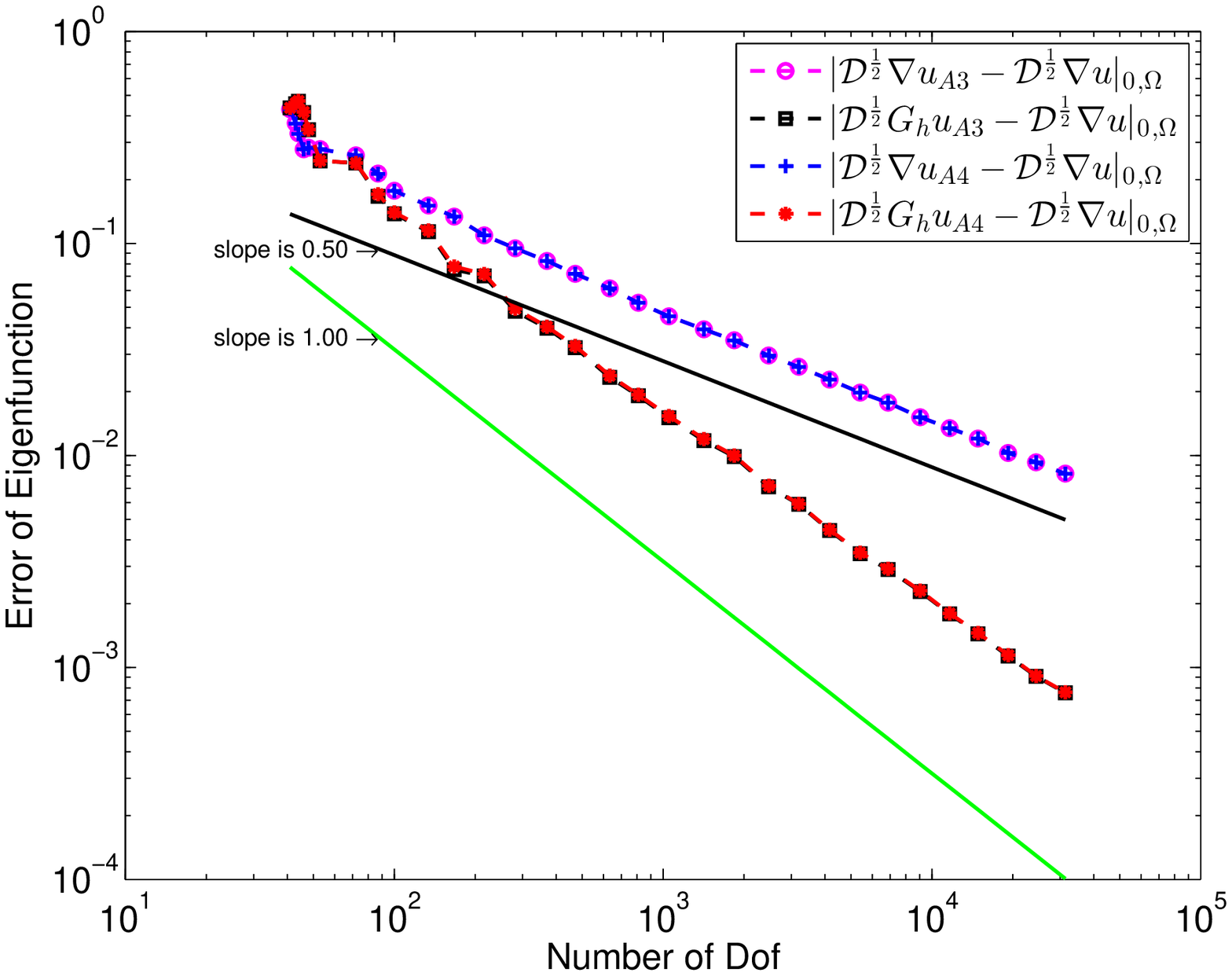}
  \caption{ Eigenfunction approximatio Error for Example 3}
  \label{fig:Electric_eigenfunction}
\end{minipage}
\end{figure}

We solve this eigenvalue problem with $\Omega = (-5, 5) \times (-5, 5)$
and zero boundary condition as in \cite{xie}. The initial mesh is shown
in Fig  \ref{fig:Electric_initial} and the adaptive mesh after 20
iterations is displayed in Fig \ref{fig:Electric_adaptive}.
The parameter $\theta $ is chosen as $0.4$.  Numerical results
are presented in Figs \ref{fig:Electric_eigenvalue} and
\ref{fig:Electric_eigenfunction}. For eigenvalue approximation,
$O(N^{-1})$ convergence rate is observed for $|\bar{\lambda}_{\text{A3}}-\lambda|$
while $O(N^{-2})$ ultraconvergence  rate is observed for $|\lambda_{\text{A3}}-\lambda|$.
For eigenfunction approximation, $\|\mathcal{D}^{\frac{1}{2}}
\nabla u_{\text{A3}} - \mathcal{D}^{\frac{1}{2}} \nabla u\|_{0, \Omega}
\approx O(N^{-0.5})$ and $\|\mathcal{D}^{\frac{1}{2}}
G_h u_{\text{A3}} - \mathcal{D}^{\frac{1}{2}} \nabla u\|_{0, \Omega}
\approx O(N^{-1})$. The numerical result of Algorithm \ref{alg:improve}
is similar.

\begin{figure}[!h]
  \centering
  \begin{minipage}[c]{0.5\textwidth}
  \centering
  \includegraphics[width=0.9\textwidth]{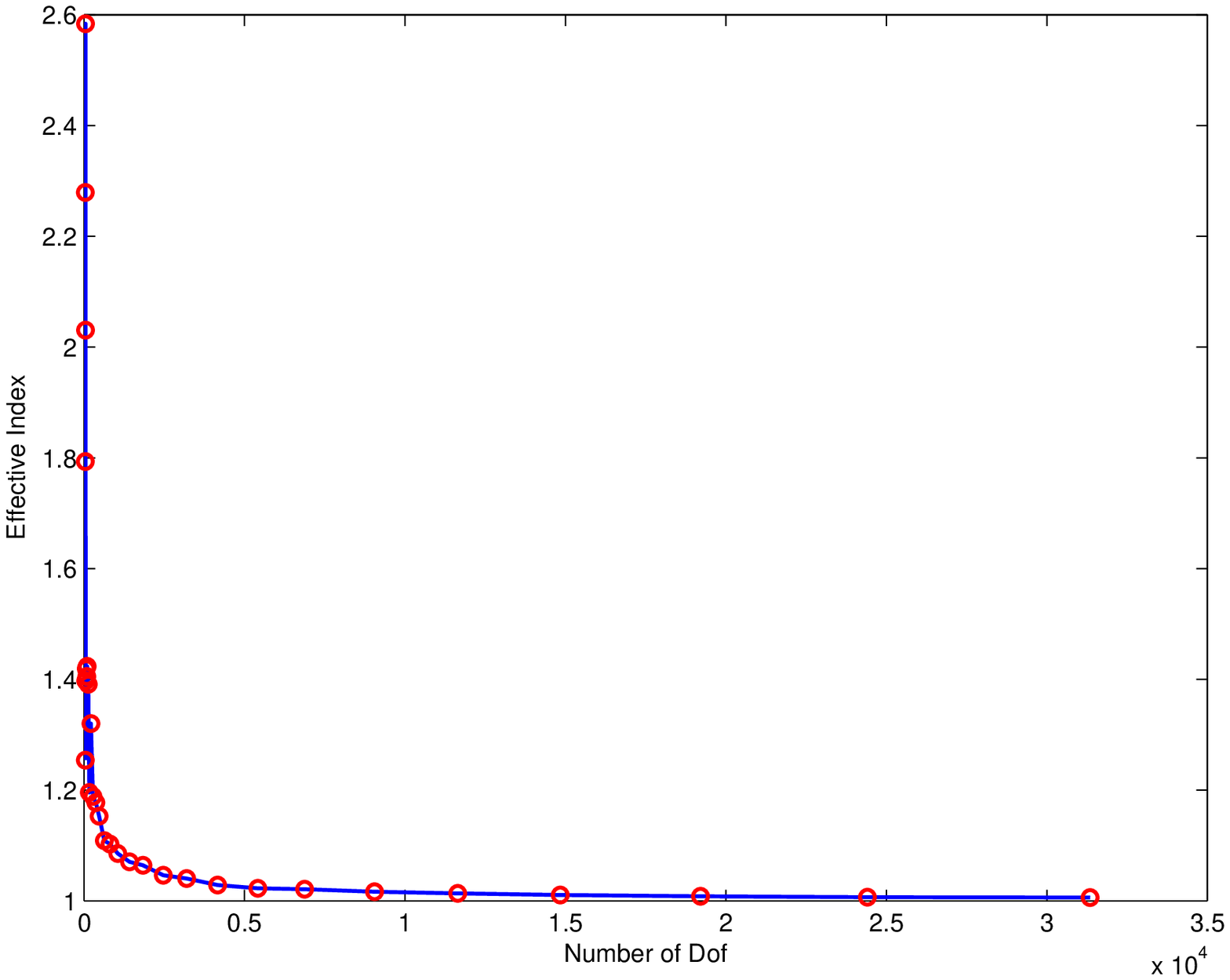}
  \caption{Effective index of Algorithm \ref{alg:adaptive} for Example 3}
\label{fig:effectindexa3}
\end{minipage}%
 \begin{minipage}[c]{0.5\textwidth}
  \centering
  \includegraphics[width=0.9\textwidth]{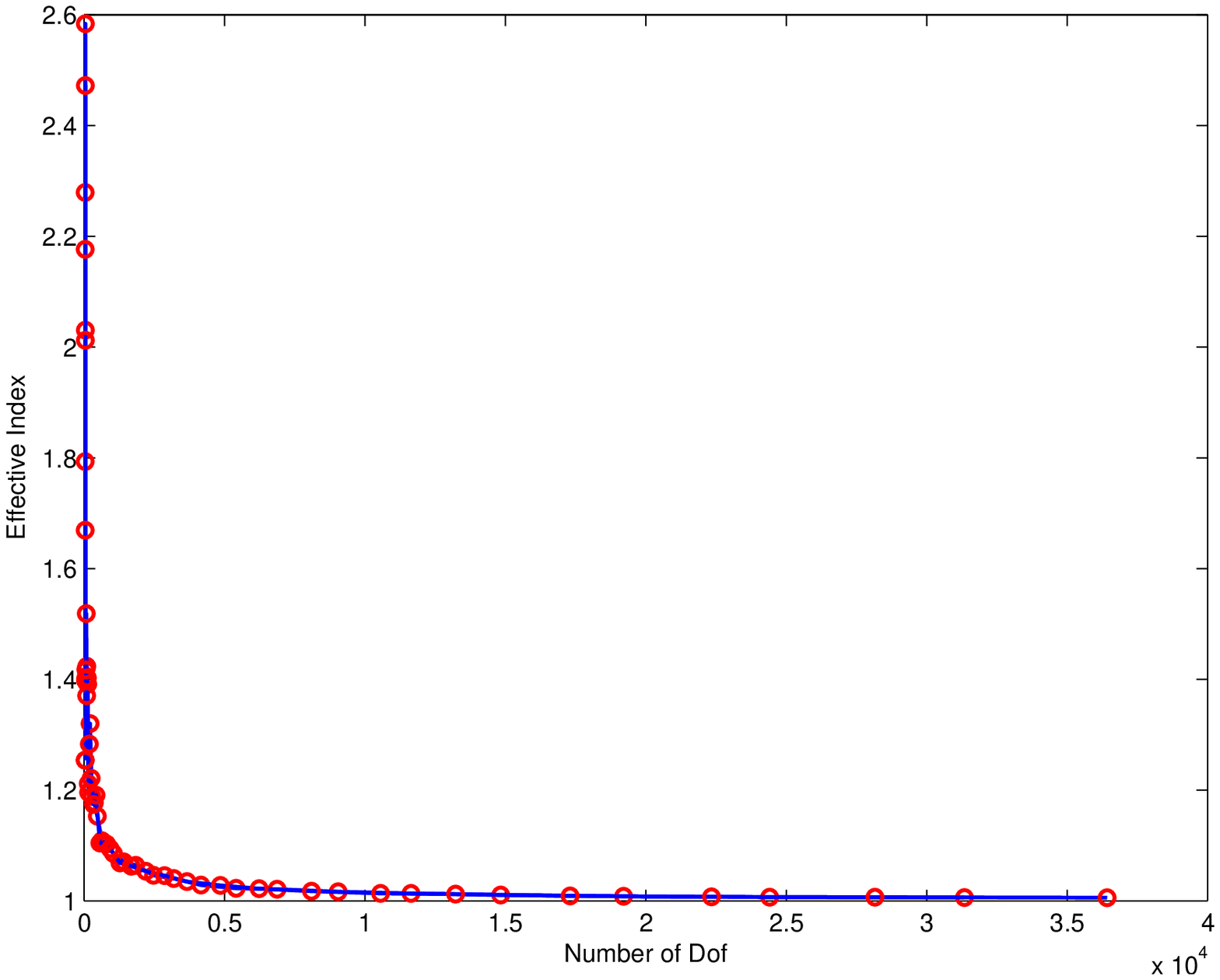}
  \caption{Effective index of Algorithm \ref{alg:improve} for Example 3}
\label{fig:effectindexa4}
\end{minipage}
\end{figure}

Figs \ref{fig:effectindexa3} and \ref{fig:effectindexa4} graph
the eigenvalue effectivity index for
the two proposed multilevel adaptive algorithms.
It also indicates that the posteriori error estimator \eqref{equ:indicator} or
\eqref{equ:globalindicator} is asymptotically exact for 
problem \eqref{equ:electric}.

\section{Conclusion}
When eigenfunctions are relatively smooth, two-space methods (using higher-order elements in the second stage)
is superior to two-grid methods (using the same element at finer grids in the second stage). They have the
comparable accuracy. However, at the last stage,
the degrees of freedom of the two-space method is much smaller than that of the two-grid method.

For linear element on structured meshes, using gradient recovery at the last stage achieves similar accuracy as the
quadratic element on the same mesh. Therefore, with much reduced cost, the gradient recovery is comparable with the two-stage method on the same mesh.

Algorithms 3 and 4 use recovery type error estimators to adapt the mesh,
and have two advantages comparing with the residual based adaptive algorithms.
1) Cost effective. In fact, the recovery based error estimator plays two roles:
one is to measure the error, and another is to enhance the eigenvalue approximation.
2) Higher accuracy. Indeed, after recovery enhancement, the approximation error is further reduced.


\end{document}